\documentclass[]{spie}  

\usepackage[superscript]{cite}

\usepackage{amssymb}
\usepackage{amsmath}
\usepackage{relsize}
\usepackage{graphicx}
\usepackage{verbatim}
\usepackage{enumitem}
\usepackage[colorlinks=true, allcolors=blue]{hyperref}
\usepackage{cite}
\usepackage{amsfonts}
\usepackage{CJK}
\usepackage{multirow}
\usepackage{diagbox}

\DeclareFontFamily{U}{mathx}{\hyphenchar\font45}
\DeclareFontShape{U}{mathx}{m}{n}{
      <5> <6> <7> <8> <9> <10>
      <10.95> <12> <14.4> <17.28> <20.74> <24.88>
      mathx10
      }{}
\DeclareSymbolFont{mathx}{U}{mathx}{m}{n}
\DeclareFontSubstitution{U}{mathx}{m}{n}
\DeclareMathAccent{\widecheck}{0}{mathx}{"71}

\newcommand{\tbf}{\textbf}
\newcommand{\eps}{\epsilon}
\newcommand{\Sch}{S}
\newcommand{\oo}{\infty}

\newcommand{\R}{\mathbb{R}}
\newcommand{\N}{\mathbb{N}}
\newcommand{\Z}{\mathbb{Z}}

\newcommand{\C}{\mathbb{C}}

\newcommand{\rar}{\rightarrow}

\newcommand{\pt}{\partial}

\newcommand{\charone}{\textbf{\large{1}}}

\newcommand{\sinc}{\text{sinc}}

\newcommand{\twopi}{2\pi}

\newcommand{\supp}{\text{supp}}

\newcommand{\Bepsd}{\mathcal{B}_{\eps,d}}

\newcommand{\epstilde}{\tilde{\eps}}
\newcommand{\Nset}{\Z_N}
\newcommand{\Lset}{\Z_L}
\newcommand{\Kset}{\Z_K}

\newcommand{\Tstep}{T_\text{step}}
\newcommand{\Ystep}{Y_\text{step}}
\newcommand{\Tshift}{T_\text{shift}}
\newcommand{\Yshift}{Y_\text{shift}}
\newcommand{\vskipp}{\vspace*{2.5mm}}

\newcommand{\im}{{\mathrm i}}
\newcommand{\e}{{\mathrm e}}
\newcommand{\dd}[1]{\, {\mathrm d}{#1}}

\newcommand{\STFT}{\mathcal{V}_g f}

\newcommand{\QSTFT}{\mathcal{V}^Q_{h}f}

\newcommand{\SSTSTFT}{\mathcal{S}^{\beta}_{g,\gamma}f}

\newcommand{\SSTSTFTepstilde}{\mathcal{S}^{\beta}_{g,\tilde{\eps}}f}

\newcommand{\SSTQSTFT}{\mathcal{S}^{Q,\beta}_{h,\gamma}f}
\newcommand{\SSTQSTFTd}{S^{Q}_{h,\gamma}f}
\newcommand{\SSTQSTFTdMAG}{S^{|Q|^2}_{h,\gamma} f}
\newcommand{\SSTQSTFTepstilde}{\mathcal{S}^{Q,\beta}_{h,\tilde{\eps}}f}

\newcommand{\IFinfoSTFT}{\Xi_{g}f}

\newcommand{\IFinfoQSTFT}{\Xi^Q_{h}f}

\newcommand{\AsetSTFT}{{\mathcal{A}_{g,\gamma}f(t)}}

\newcommand{\AsetQSTFT}{{\mathcal{A}^Q_{h,\gamma}f(t)}}

\newcommand{\AsetQSTFTd}{{A^Q_{h,\gamma}f[n]}}
\newcommand{\BsetQSTFTd}{{B^Q_{h,k}f[n]}}

\newcommand{\RTSTFT}{\tau_{g}f}

\newcommand{\STFTd}{V_gf}
\newcommand{\STFTdfp}{V^+_gf}

\newcommand{\QSTFTd}{V^Q_hf}
\newcommand{\QSTFTdhw}{V^Q_{h^w}f}
\newcommand{\QSTFTdp}{V^{Q}_{h}f^+}
\newcommand{\IFinfoSTFTd}{\widetilde{\Xi}_{g}f}
\newcommand{\IFinfoSTFTdfplus}{\widetilde{\Xi}_{g}f^+}
\newcommand{\IFinfoSTFTdfminus}{\widetilde{\Xi}_{g}f^-}
\newcommand{\IFinfoSTFTdsecond}{\widetilde{\Xi}^{(2)}_{g}f}

\newcommand{\IFinfoQSTFTd}{\widetilde{\Xi}^Q_{h}f}
\newcommand{\IFinfoQSTFTdInverse}{(\widetilde{\Xi}^Q_{h}f)^{-1}}

\newcommand{\RTSTFTd}{T_{g}f}
\newcommand{\RTSTFTdfplus}{T_{g}f^+}
\newcommand{\RTSTFTdfminus}{T_{g}f^-}
    
\newcommand{\itemsubeqn}{\hfill\refstepcounter{equation}\textup{(\theequation)}}

\newtheorem{defn}[theorem]{Definition}
\let\olddefn\defn
\renewcommand{\defn}{\olddefn\normalfont}
\newtheorem{remark}[theorem]{Remark}
\let\oldremark\remark
\renewcommand{\remark}{\oldremark\normalfont}
\newtheorem{notation}{Notation}

\let\oldnotation\notation
\renewcommand{\notation}{\oldnotation\normalfont}
\newtheorem{algorithm}{Algorithm}
\let\oldalgorithm\algorithm
\renewcommand{\algorithm}{\oldalgorithm\normalfont}
\newcommand{\citenote}[2]{%
\mbox{\cite{#1}%
\textsuperscript{,\,#2}}%
}

\title{Adaptive synchrosqueezing based on a quilted short-time Fourier transform}

\author{Alexander Berrian\footnote{\hspace*{2mm}Email: \url{aberrian@math.ucdavis.edu} -- Website: \url{http://www.math.ucdavis.edu/\textasciitilde aberrian}} \hspace*{0.5mm} and Naoki Saito\footnote{ \hspace*{1mm} Email: \url{saito@math.ucdavis.edu} -- Website: \url{http://www.math.ucdavis.edu/\textasciitilde saito}  } }
\affil{Department of Mathematics, University of California, One Shields Avenue, Davis, CA 95616, USA}

\begin{document}  
\maketitle    

\begin{abstract}
In recent years, the \emph{synchrosqueezing transform} (SST) has gained popularity as a method for the analysis of signals that can be broken down into multiple components determined by instantaneous amplitudes and phases.  One such version of SST, based on the short-time Fourier transform (STFT), enables the sharpening of instantaneous frequency (IF) information derived from the STFT, as well as the separation of amplitude-phase components corresponding to distinct IF curves.  However, this SST is limited by the time-frequency resolution of the underlying window function, and may not resolve signals exhibiting diverse time-frequency behaviors with sufficient accuracy.  In this work, we develop a framework for an SST based on a ``quilted'' short-time Fourier transform (SST-QSTFT), which allows adaptation to signal behavior in separate time-frequency regions through the use of multiple windows.  This motivates us to introduce a discrete reassignment frequency formula based on a finite difference of the phase spectrum, ensuring computational accuracy for a wider variety of windows.  We develop a theoretical framework for the SST-QSTFT in both the continuous and the discrete settings, and describe an algorithm for the automatic selection of optimal windows depending on the region of interest.  Using synthetic data, we demonstrate the superior numerical performance of SST-QSTFT relative to other SST methods in a noisy context. Finally, we apply SST-QSTFT to audio recordings of animal calls to demonstrate the potential of our method for the analysis of real bioacoustic signals.
\end{abstract}

\section{Introduction}
In the field of signal processing, one frequently seeks to model a signal with time-varying oscillatory properties as a sum of distinct \emph{amplitude-phase components} containing information on instantaneous amplitudes (IAs) and instantaneous frequencies (IFs).  In particular, one might characterize a signal $f: \R \rar \C$ using an \emph{amplitude-phase decomposition} given by $f(t) = \sum_{m=1}^M f_m(t); \,\,\, f_m(t) := A_m(t) \e^{\twopi\im \phi_m(t)},$ where the $\{A_m\}$ represent \emph{instantaneous amplitudes} (IAs), the $\{\phi_m\}$ represent \emph{instantaneous phases} (IPs), and the  $\{\phi'_m\}$ represent \emph{instantaneous frequencies} (IFs). Then, the problem is to retrieve the IFs $\phi'_m$ and amplitude-phase components (\emph{modes}) $f_m$, given that only $f$ is known.  

One may use a time-frequency representation such as the short-time Fourier transform (STFT) or continuous wavelet transform (CWT) to analyze the signal.  However, these transforms provide blurry amplitude and frequency information, thereby complicating the task of accurately determining each amplitude-phase component. Hence, often a post-processing method is used to sharpen the blurry signal information.  The tool of interest to us is the \emph{synchrosqueezing transform} (SST).  Originally introduced by Daubechies and Maes \cite{daubechies1996nonlinear} in the context of audio signal processing, the SST is a time-frequency representation that provides a sharpened picture of IAs and IFs and enables the separation and reconstruction of each separate amplitude-phase component for a certain class of signals \cite{daubechies2011synchrosqueezed, thakur2011synchrosqueezing, oberlin2014fourier}. 

The SST technique has been applied to problems in many different disciplines, and researchers have further developed the SST idea to work in an increasingly large collection of signal processing contexts.  In 2011, Daubechies, Lu, and Wu \cite{daubechies2011synchrosqueezed} gave theoretical proofs of the effectiveness of the CWT-based version of the SST and demonstrated its applicability to sets of medical data.  The same year, Thakur and Wu \cite{thakur2011synchrosqueezing} generalized the SST to a version based on STFT. Since then, there has been an uptick of research on SST, and the technique has been successfully adapted to a number of physical problems, including speaker identification from an audio signal \cite{daubechies1996nonlinear}, fault diagnosis in planetary gearboxes for wind turbines \cite{feng2015iterative}, quantifying the effect of solar radiation on a paleoclimate change on Earth \cite{thakur2013synchrosqueezing}, and extracting heart-rate variability (HRV) from an ECG signal \cite{daubechies2011synchrosqueezed}.

However, researchers continue to develop techniques for cases where the original SST techniques are insufficient.  Often, this insufficiency is due to the resolution constraints of the underlying time-frequency representation (STFT or CWT).    One area of recent work can be categorized under the umbrella of \emph{non-stationary Gabor transforms} (NSGT)\cite{balazs2011theory}, where one analyzes the signal using a kernel function with changing time-frequency resolution over the time-frequency plane.  This allows for generalizations of the STFT and CWT that allow the user to adapt the time-frequency resolution to the signal content. Balazs et al.\cite{balazs2011theory}\ developed the idea of NSGT in a frame theory context, and D\"{o}rfler further generalized this idea to the notion of \emph{quilted Gabor frames} \cite{dorfler2011quilted}. Her work refers to the time-frequency plane as a ``quilt'' with different ``patches'' corresponding to regions where a signal exhibits different time-frequency behavior.

In the context of SST, recent developments under the umbrella of NSGT include generalizations of SST where different windows may be used to adapt to the signal at different times \cite{guo2014new}, or at different frequencies \cite{holighaus2016reassignment}.  In a previous work\cite{berrian2015quilted}, we further generalized this idea by introducing an SST based on a \emph{quilted STFT} (QSTFT), where the window is allowed to change depending on the time-frequency region of interest.  More recently, Sheu et al.\ derived an algorithm for SST with the automatic selection of optimal windows in different time-frequency regions\cite{sheu2015entropy}, but did not provide a method for signal reconstruction in the context of joint time-frequency window variation.  In this work, we provide theoretical results for the effectiveness of the SST-QSTFT in both the continuous and the discrete context, and we derive a slightly different algorithm for automatic window selection.  Furthermore, we explore the usage of chirped windows in depth and provide a methodology for signal reconstruction when the selected windows vary in both time and frequency. 

We give an outline of this paper as follows. In Sec.\ 2, we review the relevant background material on the STFT and the corresponding synchrosqueezing transform.  We then introduce the QSTFT and SST-QSTFT in the continuous setting (Sec.\ 3) and in the discrete setting (Sec.\ 4), and we provide theorems that demonstrate the effectiveness of the SST-QSTFT for accurate IF detection and mode reconstruction.  Numerical results and applications follow in Sec.\ 5.  We conclude with a summary of our results in Sec.\ 6.

\section{Background}

\subsection{Admissible function class}
In this analysis, we restrict our study to a class of signals satisfying certain theoretical properties.
\begin{defn}
A signal $f$ is said to be in the \emph{weakly modulated IA and IF signal class} $\Bepsd$ \cite{oberlin2014fourier, thakur2011synchrosqueezing, daubechies2011synchrosqueezed, thakur2013synchrosqueezing} if for some $M \in \N$ we can write 
\begin{equation}\label{amp-phase-decomp}
f(t) = \sum_{m=1}^M f_m(t); \,\,\, f_m(t) := A_m(t) \e^{\twopi\im \phi_m(t)};
\end{equation}
and if there exist $\eps, d > 0$ such that for each $m \in \{1,\dots,M\}$,
\begin{itemize}
\item $A_m$ and $\phi'_m$ are \emph{bounded and sufficiently smooth}: $A_m \in C^1\cap L^\oo$, $\phi_m \in C^2$, $\phi'_m \in L^\oo$, $\inf\limits_{t\in\R} A_m(t) > 0$ and $\inf\limits_{t\in\R} \phi'_m(t) > 0$;
\item the IA $A_m$ and IF $\phi'_m$ are \emph{slowly-varying}: $\|A_m'\|_\oo \leq \eps$ and $\|\phi_m''\|_\oo \leq \eps$; 
\item $\phi'_m$ is \emph{well-separated} from the other IFs: $\phi'_m(t) - \phi'_{m-1}(t) > d$ for each $t \in \R$, provided that $2 \leq m \leq M$. 
\end{itemize}
\end{defn}
\noindent The assumptions above are necessary to ensure the accurate isolation of the IFs $\{\phi'_m\}$ and reconstruction of the modes $\{f_m\}$ in the theory that follows.
\begin{remark}
One models real-valued signals by the real part in (\ref{amp-phase-decomp}), yielding $f_\text{real}(t) = \displaystyle \sum_{m=1}^M A_m(t) \cos(\twopi\phi_m(t))$. 
\end{remark}

\subsection{Fourier transforms}
The standard building block for determining the spectral information of a signal is the \emph{Fourier transform}, which comes in different forms and can be defined in many different ways. In this paper, we define the \emph{continuous Fourier transform} of a continuous-time signal $f \in L^2(\R) \cap L^1(\R)$ by $\hat{f}(\xi) := \displaystyle \int_\R f(x) \e^{-\twopi\im x\xi}\dd{x}$ for each $\xi \in \R$.  Also, for the purposes of this paper, we define the \emph{semi-discrete Fourier transform} of the periodic sequence $g$ by $\hat{g}(u) := \displaystyle \sum_{\ell=0}^{L-1} g[\ell] \e^{-\twopi\im u \ell}$ for the continuous argument $u \in \R$.\vskipp
\begin{remark}
It is a slight abuse of notation to use the same hat symbol for the continuous Fourier transform $\hat{f}$ of a continuously-time function $f$ and for the semi-discrete Fourier transform $\hat{g}$ of a (discrete) periodic sequence $g$.  However, the separate meaning of these two transforms is clear based on whether the transform input is continuously or discretely valued.  
Moreover, in the rest of this paper, we will only consider the Fourier transform of continuously-time window functions and discrete periodic window sequences as defined in Sec.\ \ref{section-stft}, and we will not use the hat symbol in ambiguous cases where (for instance) a discrete periodic window sequence is explicitly said to be derived from discretizing a continuous-time window function.
\end{remark} 

\subsection{Short-time Fourier transform (STFT)}\label{section-stft}
In order to analyze the local frequency content at different times of the signal $f$, we first consider the \emph{continuous short-time Fourier transform} of $f$ with respect to a window function $g \in L^2(\R)$ centered at $0$, which we define \cite{muller2015fundamentals} by
\begin{align*}
\STFT(t,\xi) &:= \int_\R f(x) \overline{g(x-t)} \e^{-\twopi\im\xi(x-t)}\dd{x}.
\end{align*}
Here, the effect of the window $g$ is to essentially truncate the signal $f$ around the time of interest $t$ in a smooth fashion, in order to enable the user to examine the local spectral information of the signal near the time $t$.  If the window $g$ is compactly supported, then we say $g$ is \emph{time-limited}.  If $\hat{g}$ is compactly supported, then we say $g$ is \emph{band-limited}.  Note that a function $h$ is said to be \emph{compactly supported} if the set $\supp(h) := \{ x \in \R \ : \ h(x) \neq 0\}$ satisfies $\supp(h) \subseteq [-R, R]$ for some finite $R > 0$.

In the context of discrete implementation, we assume that $f$ is sampled at sampling rate $f_s$, denoting $\Delta t := 1/f_s$.  Here and in the rest of this paper, we abuse notation slightly by denoting $f[\ell] := f(\ell\Delta t)$, $f_m[\ell] := f_m(\ell\Delta t)$, $A_m[\ell] := A_m(\ell\Delta t)$, $A_m'[\ell] := A_m'(\ell\Delta t)$, and similarly for $\phi_m$ and all its derivatives.\footnote{In this paper, we adhere to the convention of using square brackets for functions with discrete inputs, and parentheses for functions with continuous inputs.}   We also further assume that $f$ is strictly zero outside a compact region having left boundary at index $\ell=0$, as it is practical to assume for many real-world applications where we are only interested in a finite-length signal.

We may now define the \emph{discrete short-time Fourier transform} of $f$ for a discrete, time-limited window function $g$ sampled at $L$ points, with $g[\ell] := g_\ell$ for $\ell = 0, \dots, L-1$, and with \emph{hop size} $H \in \N$, via\cite{muller2015fundamentals}
\begin{equation}\label{discrete-STFT}
\STFTd[n,k] := \sum_{\ell=0}^{L-1} f[\ell + nH] \overline{g[\ell]} \e^{-\twopi\im k\ell/L},
\end{equation}
for \emph{frequency bins} $k \in \{0, \dots, L-1\}$, and \emph{frames} $n \in \{0, \dots, N-1\}$ where $N \in \N$ is large enough so that $f[\ell] = 0$ for all $\ell \geq (L-1) + (N-1)H$. For signals with many samples, it is often practical to set the hop size $H$ larger than $1$, yielding a subsampling of the STFT.  In such a context, an \emph{overlap-add formula} can provide a perfect reconstruction of the original signal, provided that the window satisfies a certain formula for the given hop size \cite{SASPWEB2011}\footnote{See \url{http://ccrma.stanford.edu/~jos/sasp/Choice_Hop_Size.html} for the specific formula, and which windows satisfy it.}.

\begin{figure}
\begin{center}
\includegraphics[height=8cm,width=18cm]{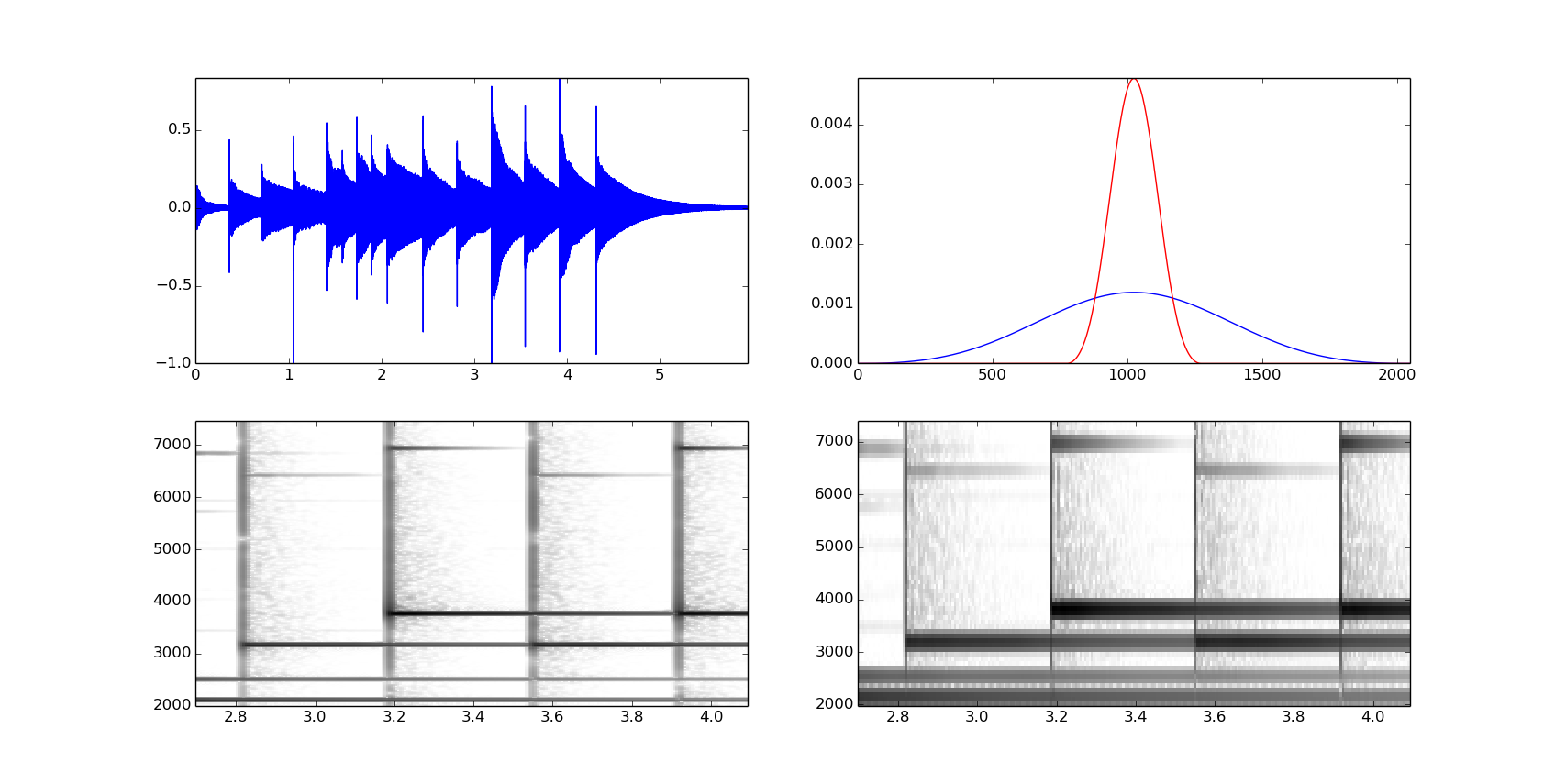}
\end{center}\vspace*{-3mm}
\caption{\label{figure-two-stfts} \tbf{Top left:} Glockenspiel signal $f$, sampled at $f_s := 44100$ Hz. \tbf{Top right:} Window functions used to analyze $f$. \emph{Short window (in blue):}  Blackman window $g^1$ of length 2000 samples. \emph{Tall window (in red):}  Blackman window $g^2$ of length 500 samples. $g^1$ and $g^2$ are normalized to have unit $\ell^1$-norm and zero-padded to 2048 and 512 samples respectively. \tbf{Bottom left:} $\log_2(1 + |V_{g^1}f[n,k]|^2)$, for times $2.7$ to $4.1$ seconds and frequencies $2000$ to $7500$ Hz. \tbf{Bottom right:} $\log_2(1 + |V_{g^2}f[n,k]|^2)$, over the same frequency range.  $V_{g^1}f$ and $V_{g^2}f$ are both computed with hop size $H = 250$ samples. In the STFT plots, the vertical content corresponds to transient events (note onsets) and the horizontal content corresponds to tonal events (enduring notes). With $g^1$, the frequency resolution is better, while the time resolution is worse, leading to better visualization of note pitches but worse visualization of note onset times. With $g^2$, the time resolution is better, leading to sharp onset resolution, but the frequency resolution is worse, leading to obscured pitch information.}
\end{figure}

The STFT allows one to visualize the IAs and IFs of the signal $f$, as well as the presence of transient and noise content spanning over many frequencies.  Due to the Fourier uncertainty principle\citenote{folland1992fourier}{Sec.\ 7.3}, there is a tradeoff between the time precision and the frequency precision of the STFT.  That is, a window $g$ that enables the STFT to render more precise frequency information will yield less precise time information, and vice-versa.  This point is illustrated using STFTs with respect to two different windows in Figure \ref{figure-two-stfts}.  Thus, the Fourier uncertainty principle limits the robustness of the STFT, and so various methods have been designed to sharpen the STFT information and approximately retrieve the IAs and IFs of a signal.  In the next subsection, we describe the \emph{STFT-based synchrosqueezing transform}, which enables the sharpening of the frequency information derived from the STFT. 

\subsection{Synchrosqueezing transform based on STFT (SST-STFT)}

\subsubsection{Continuous theory}
The idea of the continuous STFT-based synchrosqueezing transform \cite{thakur2011synchrosqueezing} is to calculate a \emph{reassignment frequency} $\IFinfoSTFT$ for each STFT coefficient with magnitude above a certain tolerance. This reassignment frequency is an estimate of the actual frequency location of the nearest IF.  Then, for each time-frequency point $(t,\xi)$, the synchrosqueezing transform finds all other frequencies $\eta$ where the STFT coefficient $\STFT(t,\eta)$ at this time has reassignment frequency $\IFinfoSTFT(t,\eta) = \xi$, and approximately sums up all such STFT coefficients to yield the SST-STFT coefficient $\SSTSTFT(t,\xi)$.  (We will explain the meaning of ``approximately'' summing shortly.) The effect is to ``squeeze'' the blurred-out STFT visualization in the frequency direction, leading to a better-concentrated time-frequency representation along the frequency axis.  The SST-STFT is defined as follows.
\begin{defn}
The \emph{continuous STFT-based synchrosqueezing transform} (SST-STFT) \cite{thakur2011synchrosqueezing} with tolerance $\gamma \geq 0$ and limiting parameter $\beta > 0$ is given by
\begin{align*}
\SSTSTFT(t,\xi) := \int_{\AsetSTFT}\STFT(t,\eta) \dfrac{1}{\beta}b\left(\dfrac{\xi - \IFinfoSTFT(t,\eta)}{\beta}\right) \dd{\eta},
\end{align*}
where $b \in C^{\oo}_c(\R)$ is
a ``bump function'' satisfying $\hat{b}(0) = 1$,
$\AsetSTFT := \{\eta \in \R^+ : |\STFT(t,\eta)| > \gamma\}$, and $\IFinfoSTFT(t,\eta) := \dfrac{\pt_t\STFT(t,\eta)}{\twopi\im \STFT(t,\eta)}$ is the \emph{STFT-based reassignment frequency}.
\end{defn}\vskipp

\begin{remark}
As $\beta \searrow 0$, the term $\frac{1}{\beta}b\left(\frac{\xi - \IFinfoSTFT(t,\eta)}{\beta}\right)$ converges in the distributional sense to $\delta\left(\xi - \IFinfoSTFT(t,\eta)\right)$, where $\delta$ denotes the Dirac delta.  This is the meaning of the ``approximate summation'' described earlier.
\end{remark}\vskipp

\begin{remark}
The choice of $\IFinfoSTFT$ is motivated by the fact that for a constant chirp $f(t) := A \e^{\twopi \im c t}$, with $A, c > 0$, we have exactly $\IFinfoSTFT(t,\eta) = c$ for all $(t,\eta)$\cite{daubechies2011synchrosqueezed, daubechies1996nonlinear}.
\end{remark}\vskipp

\begin{remark}
Here, $C^\oo_c(\R)$ denotes the class of functions where derivatives of all orders exist and are continuous (a property which we call \emph{$C^\oo$-smoothness}) and where each function is compactly supported.
\end{remark}\vskipp

\noindent The following theoretical result ensures the accuracy of the SST-STFT:\vskipp
\begin{theorem}\label{theorem-sst-stft} (Thakur \& Wu\cite{thakur2011synchrosqueezing}; Oberlin, Meignen \& Perrier \cite{oberlin2014fourier}) 
Let $\eps > 0$, $\nu \in (0,1/2)$, $\epstilde := \eps^{\nu}$, $d > 0$. Suppose that $\displaystyle f = \sum_{m=1}^M f_m \in \Bepsd$. Let $g \in \Sch(\R)$, where $\Sch(\R)$ denotes the \emph{Schwartz class} of $C^\oo$-smooth, rapidly-decaying functions. Assume that $g$ is real-valued and satisfies $\text{\supp}(\hat{g}) \subseteq [-d/2, d/2]$. Then, if $\eps$ is sufficiently small:
\begin{itemize}
\item \tbf{(Concentration of STFT around IF curves)} $|\STFT(t,\xi)| > \epstilde$ only when there is an $m \in \{1, \dots, M\}$ such that $(t,\xi) \in \mathcal{Z}_m := \{ (t, \xi) \in \R \times \R^+:\ |\phi'_m(t) - \xi| < d/2 \}$.
\item \tbf{(Closeness of reassignment frequency $\IFinfoSTFT$ to nearby IF)} For all $m \in \{1, \dots, M\}$ and all $(t,\xi) \in \mathcal{Z}_m$ such that $|\STFT(t,\xi)| > \epstilde$, we have $|\IFinfoSTFT(t,\xi) - \phi'_m(t)| \leq \epstilde$.  
\item \tbf{(Accuracy of reconstruction)} For each $m \in \{1, \dots, M\}$ there is a constant $C_m > 0$ such that for all $t \in \R$, 
\begin{align}\label{inequality-sst-stft-reconstruction}
\left| \lim_{\beta \rar 0^+} \left(\int\limits_{ \{\xi \ : \ |\xi - \phi'_m(t)| < \epstilde \}} \dfrac{1}{g(0)} \cdot \SSTSTFTepstilde(t,\xi) \dd{\xi} \right) - f_m(t)\right| &\leq C_m \epstilde.
\end{align}
\end{itemize}
\end{theorem}

\begin{remark}
The assumption that $g$ is real-valued in Theorem \ref{theorem-sst-stft} is not essential.  In Sec.\ \ref{section-theorem-continuous}, we will state a generalization of Theorem \ref{theorem-sst-stft} (Theorem \ref{theorem-continuous}), from which the statement of Theorem \ref{theorem-sst-stft} will follow for complex-valued $g$, provided that one normalizes by $\overline{g(0)}$ instead of $g(0)$ in the reconstruction formula given in (\ref{inequality-sst-stft-reconstruction}).
\end{remark}\vskipp

\subsubsection{Discrete theory}\label{section-sst-stft-discrete-theory}
The discrete theory of the SST-STFT was partially addressed in the work of Thakur and Wu\cite{thakur2011synchrosqueezing}, resulting in a theorem concerning the concentration of the STFT around the IF curves and the closeness of the reassignment frequency to the nearby IF.  The authors additionally showed robustness to noise and accuracy in the setting of nonuniform sampling.  However, the theory prescribes a reassignment frequency formula that may lead to aliasing, for the following reason. Thakur and Wu, in computing the analogue of $\pt_t V_gf$ for the discrete case, passed the derivative $\pt_t$ through the integral and computed
\begin{align*}
\pt_t \STFT(t,\xi) &= \int_\R f(x) (\twopi \im \xi g(x-t) - g'(x-t)) \e^{-\twopi\im\xi(x-t)} \dd{x},
\end{align*}
which enables one to derive an analogous discrete reassignment frequency formula to replace $\pt_t \STFT$.  However, if $g'$ is nonzero at the boundaries of the support of $g$ (i.e., if $g'$ does not tail off to zero when $g$ tails off to zero), then an aliased discrete Fourier transform (DFT) will result, since the signal is assumed to be periodic.  Moreover, the explicit form of $g'$ may not be available.  Indeed, one may wish to use a discrete-time window sequence that is not simply the discretized version of a continuous-time window function, in which case it does not make sense to talk about the continuous-time derivative $g'$.  Examples of such window sequences include the minimal-latency windows used by Su and Wu for real-time SST\cite{su2016minimum}.  These considerations lead us to suggest a reassignment frequency formula based on finite differencing of the phase spectrum in Sec.\ \ref{section-alternative-reassignment-frequency}.  We also address the issue of reconstruction in Sec.\ \ref{section-discrete-time-limited-reconstruction}.


\section{SST-QSTFT: continuous setting}
Since the STFT and SST-STFT only permit a single window choice, these transforms are limited in their capability to adapt to signals whose behavior changes depending on the time-frequency region. Due to the Fourier uncertainty principle, the window cannot have both good time resolution and good frequency resolution.  Moreover, while the SST-STFT is designed to sharpen frequency resolution, it does not improve  time resolution, and the reassignment frequency is generally less accurate for signals with fast-varying instantaneous frequencies.  Hence we consider the notion of \emph{adaptive} time-frequency transforms defined from a family of window functions, where different windows are used for time-frequency regions containing different phenomena.  We call such a family a \emph{quilted window family}, borrowing the notion of a time-frequency ``quilt'' from D\"{o}rfler, who coined the term in the context of her work on \emph{quilted Gabor frames}\cite{dorfler2011quilted}.  Here, different ``patches'' of the quilt refer to time-frequency regions with different signal behavior.   We define a modified version of the STFT for the case of a quilted window family.  Next, we define the SST based on this ``quilted'' STFT.  We then define the notion of an \emph{adaptive quilted window family} for a signal $f$ of the class $\Bepsd$.  Then we state a theorem for the theoretical accuracy of the reassignment frequency and mode reconstruction for the SST based on quilted STFT, in the case when an adaptive quilted window family is used.  Due to the page limitation, we leave the proof of this theorem for our future work.

\subsection{Continuous quilted window families}
\begin{defn}
We define a \emph{continuous quilted window family} to be a collection $\{h_{t,\xi}\}_{(t,\xi) \in \R \times \R^+}$ where for each $(t,\xi) \in \R\times\R^+$, $h_{t,\xi} \in L^2(\R)$ is a window function. Hence, we associate a window function to each time-frequency point.
\end{defn}\vskipp

\noindent As a simple example of a continuous quilted window family, we consider first a collection consisting of two Blackman windows $\{g^1, g^2\}$, where $g^1$ is a wide Blackman window and $g^2$ is a narrow Blackman window.  The wider window enables sharp resolution along the frequency axis, while the narrower window provides sharp resolution along the time axis.  In Figure \ref{figure-two-stfts}, we showed an example of such a collection, together with the different magnitude spectra that would be produced when computing the STFT with respect to each window.

Strictly speaking, the collection $\{g^1, g^2\}$ by itself does not define a continuous quilted window family, since we have not yet associated a window to each time-frequency point of $\R\times\R^+$.  In general, one should associate time-frequency points to window functions in an automatic fashion based on the signal content, and we describe algorithms in Sec.\ \ref{section-numerical-implementation-windows} for this purpose.  Heuristically, one may associate a narrow window to time-frequency regions containing either transient events resembling delta spikes or onsets of amplitude-phase components, in order to better resolve the onset times of these events.  Similarly, one may associate a wide window to regions containing amplitude-phase components enduring over a long time period, in order to better resolve the instantaneous frequency information.  Figure \ref{figure-two-stfts} demonstrated the capacity of the collection $\{g^1, g^2\}$ to sharply resolve these different categories of signal content for a glockenspiel signal.

The usage of a more general, larger window family consisting of several dilations of a single Gaussian window to compute SST was explored in depth by Sheu et al\cite{sheu2015entropy}.  The dilation of the Gaussian window corresponds to varying its effective bandwidth, enabling for sharp representation of a greater variety of time-frequency events.  As another example, we consider a \emph{chirped window family} of the form $\{g^\sigma\}_{\sigma \in \mathfrak{S}}$, where
\begin{align*}
g^\sigma(t) &:= g(t) \e^{\twopi \im \sigma t^2/2}
\end{align*}
for each $\sigma$ in a finite set of real numbers $\mathfrak{S}$, with $g$ a fixed window function. Here, the parameter $\sigma$ is called the \emph{chirp parameter} or \emph{chirp rate}.  The STFT $V_{g^\sigma}$ will then sharply concentrate around instantaneous frequency curve segments that can be closely approximated by $\phi'(t) \approx \sigma t + c$ for some $c \in \R$.  Families of chirped windows provide an alternative to varying the window width. Moreover, by determining which parameter $\sigma^* \in \mathfrak{S}$ yields the STFT $V_{g^{\sigma^*}}$ that best concentrates the signal content in a small region, one may directly infer that the signal content in that region is better linearized by a component of the form $\sigma^* t + c$ than $\sigma t + c$ for any other $\sigma \in \mathfrak{S}$.  Thus, chirped windows allow us an immediate estimate of the signal information that does not follow directly from a simple window dilation.  We will further explore the use of chirped window families in Sec.\ \ref{section-numerical-implementation}.

Using the definition of continuous quilted window family, one may consider a tiling of the time-frequency plane where each tile $\mathcal{T} \subseteq \R^2$ has a window $h^\mathcal{T}$ associated to it; i.e., $h_{t,\xi} = h^\mathcal{T}$ for each $(t,\xi) \in \mathcal{T}$.  Thus, it is possible to consider an optimization problem where the window $h^\mathcal{T}$ is chosen to adapt to the signal behavior in $\mathcal{T}$.  Indeed, such adaptivity will become necessary when considering the accuracy of an SST in the quilted window context, as we discuss in Secs. \ref{section-adaptive-quilted-window-families} and \ref{section-theorem-continuous}. 

\subsection{Continuous quilted short-time Fourier transform (QSTFT)}
\begin{defn}
Suppose $\{h_{t,\xi}\}_{(t,\xi)\in\R \times \R^+}$ is a continuous quilted window family, and define the function $h$ by $h(x,t,\xi) := h_{t,\xi}(x)$ for each $t,x \in \R$ and $\xi \in \R^+$. Then we define the \emph{quilted short-time Fourier transform}\cite{berrian2015quilted} of the signal $f$ with respect to $h$ by
\begin{align}\label{definition-QSTFT}
\QSTFT(t,\xi) &:= \int_\R f(x) \overline{ h_{t,\xi}(x-t)} \e^{-\twopi \im \xi (x-t)} \dd{x}.
\end{align}
\end{defn}

\subsection{Continuous QSTFT-based SST (SST-QSTFT)}
\begin{defn}
 We define the \emph{continuous QSTFT-based synchrosqueezing transform}\cite{berrian2015quilted} (continuous SST-QSTFT)
of a signal $f$, with respect to the function $h$ defining the continuous quilted window family $\{h_{t,\xi}\}_{(t,\xi) \in \R \times \R^+}$, and with tolerance $\gamma\geq 0$ and limiting parameter $\beta>0$, as follows:
\begin{displaymath}
\SSTQSTFT(t,\xi) := \int_{\AsetQSTFT}\QSTFT(t,\eta)\dfrac{1}{\beta}b\left(\dfrac{\xi - \IFinfoQSTFT(t,\eta)}{\beta}\right) \dd{\eta},
\end{displaymath}
with $b$ as before, $\AsetQSTFT := \{\eta \in \R^+ : |\QSTFT(t,\eta)| > \gamma\}$, and where $\IFinfoQSTFT(t,\xi) := \dfrac{\pt_t\QSTFT(t,\xi)}{\twopi\im \QSTFT(t,\xi)}$ is the \emph{QSTFT-based reassignment frequency.}
\end{defn}

\subsection{Adaptive continuous quilted window function families}\label{section-adaptive-quilted-window-families}
In order to ensure the accuracy of reassignment frequency and mode reconstruction for the SST-QSTFT, one must make some assumptions on the quilted window family.  In particular, one must assume that the family is \emph{adaptive} to the signal in the following manner.\vskipp

\begin{defn}
Suppose that $f \in \Bepsd$. We say that the continuous quilted window family $\{h_{t,\xi}\}_{(t,\xi)\in\R \times \R^+}$ is of class $\mathcal{W}^Q_{d,\eps,f}$ if the following conditions hold:
\begin{subequations}
\begin{itemize}
\item \tbf{Smoothness, non-triviality, and band-limitation:} For each $(t,\xi) \in \R \times \R^+$ we have $h_{t,\xi} \in \Sch(\R)$, $h_{t,\xi}(0) \neq 0$, and $\text{\supp}\left(\widehat{h_{t,\xi}}\right) \subset [-d/2, d/2]$.
\item \tbf{Window choice remains constant in the frequency band around an IF value:} For each $t \in \R$ and $m \in \{1, \dots, M\}$, there exists a single window function $g_{t,m}$ such that $h_{t,\xi} \equiv g_{t,m}$ for all $\xi$ in the frequency band $\{\xi : |\phi'_m(t) - \xi| < d/2\}$.
\item \tbf{Integration bounds:} For each $p \in \{0, 1, 2\}$, there exists $I_p \in \R^+$ such that \\ $\displaystyle \sup_{(t,\xi) \in \R \times \R^+} \int_\R |u|^p \left|h_{t,\xi}(u)\right| \dd{u} \leq I_p$. \itemsubeqn \label{definition-Ip}
\item \tbf{The window family does not change too quickly over time:} Defining $h(x;t,\xi) := h_{t,\xi}(x)$, we have $\displaystyle \sup_{(t,\xi)\in\R\times\R^+} \int_\R |\pt_t h(u;t,\xi)| \dd{u} \leq \eps$ and $\displaystyle \sup_{(t,\xi)\in\R\times\R^+} \int_\R |u| \left|\pt_t h(u;t,\xi)\right| \dd{u} \leq J_1$ for some $J_1 \in \R^+$. \itemsubeqn \label{definition-Jp}
\end{itemize}
\end{subequations}
Moreover, we call $\mathcal{W}^Q_{d,\eps,f}$ the \emph{class of $f$-adaptive continuous quilted window families}.
\end{defn}\vskipp

\begin{figure}
  \begin{minipage}[c]{0.5\textwidth}
    \includegraphics[width=\textwidth, height=0.8\textwidth]{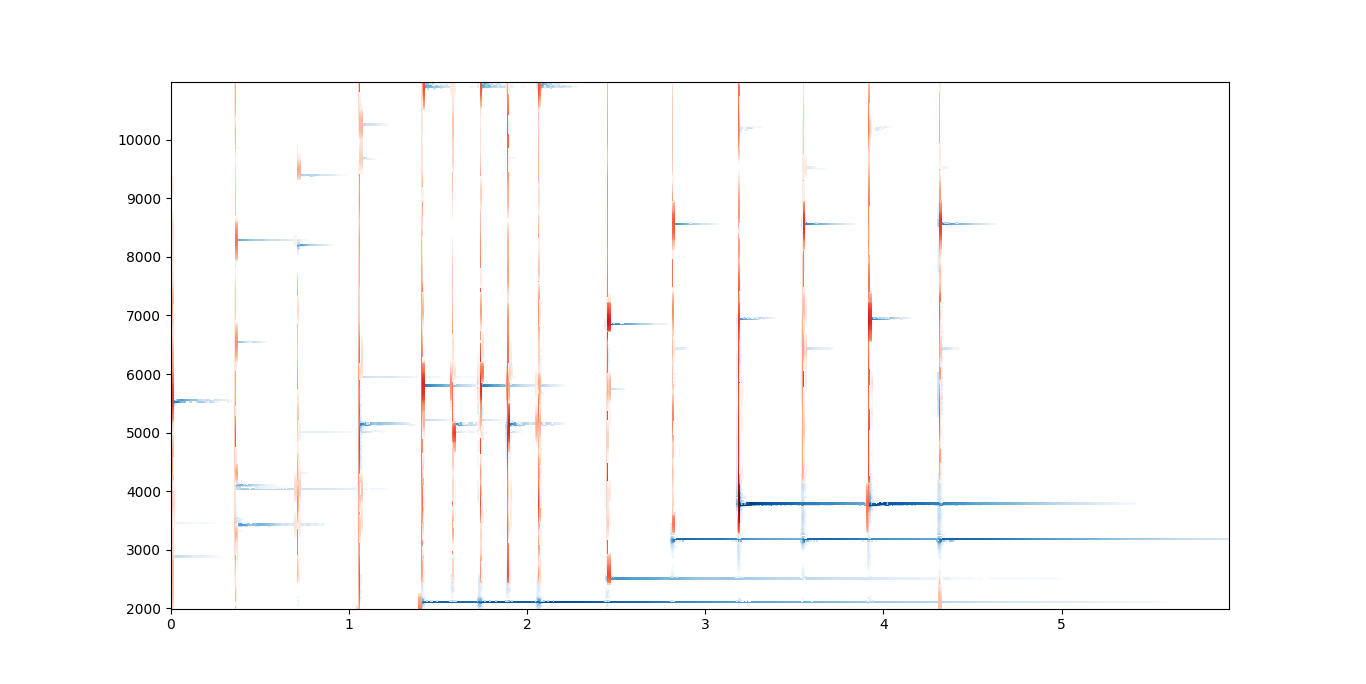}
  \end{minipage}\hfill
  \begin{minipage}[c]{0.5\textwidth}
    \caption{\label{figure-sst-qstft} SST-QSTFT of glockenspiel signal with respect to $\{\tilde{g}^1, \tilde{g}^2\}$ and with tolerance $\gamma=0$, using Algorithm \ref{algorithm-2} for automatic adaptive window selection (see Sec.\ \ref{section-numerical-implementation-windows}), with SST applied only where $\tilde{g}^1$ is used, and frequency range between $2000$ and $11000$ Hz.  $\tilde{g}^1$ and $\tilde{g}^2$ are normalized to have unit $\ell^1$-norm. Colors denote the window chosen in each time-frequency region (blue for $\tilde{g}^1$ and red for $\tilde{g}^2$). For improved energy concentration, instead of reassigning the original QSTFT coefficients, we reassign their squared magnitudes (See Secs.\ \ref{section-discrete-sst-qstft} and \ref{section-discrete-time-limited-reconstruction}). We also apply a threshold to the SST-QSTFT, to give a cleaner visualization in low-energy regions where the automatic window choice algorithm becomes unreliable. The hop size is $H = 128$ samples, and the parameters used for Algorithm \ref{algorithm-2} are $A^* = B^* = 5$, $\alpha = 0.25$, $\Tshift = \Yshift = 3$, $\Tstep = \Ystep = 1$, with $P_T = P_Y = \{-1, 1\}$.}
 \end{minipage}
\end{figure}

\begin{remark}
Of course, one does not generally know the information $\phi'_m$ in advance.  Hence, an algorithm must be prescribed to automatically select the windows in a signal-adaptive manner. We describe algorithms for this purpose in Sec.\ \ref{section-numerical-implementation-windows}.  Figure \ref{figure-sst-qstft} demonstrates the usage of an automatic adaptive window selection algorithm to compute the SST-QSTFT, using a collection of Blackman windows $\{\tilde{g}^1, \tilde{g}^2\}$, where $\tilde{g}^1$ is of length 2048 samples and $\tilde{g}^2$ is of length 512 samples.  In this figure, transient note actions are well-concentrated, owing to the better time resolution of $\tilde{g}^2$ in transient note regions. 
\end{remark}

\subsection{Theorem}\label{section-theorem-continuous}
We state a theoretical result analogous to that of Theorem \ref{theorem-sst-stft} for the continuous SST-QSTFT here:
\begin{theorem}\label{theorem-continuous}
Let $\eps > 0$, $\nu \in (0,1/2)$, $\epstilde := \eps^{\nu}$, $d > 0$. Suppose that $\displaystyle f = \sum_{m=1}^M f_m \in \Bepsd$. Assume that $\{h_{t,\xi}\}_{(t,\xi)\in\R \times\R^+}$ is of the class $\mathcal{W}^Q_{d,\eps,f}$. Then, if $\eps$ is sufficiently small we have:
\begin{itemize}
\item \tbf{(Concentration of QSTFT around IF curves)} $|\QSTFT(t,\xi)| > \epstilde$ only when there is an $m \in \{1, \dots, M\}$ such that $(t,\xi) \in \mathcal{Z}_m := \{ (t, \xi) \in \R \times \R^+:\ |\phi'_m(t) - \xi| < d/2 \}$.
\item \tbf{(Closeness of reassignment frequency $\IFinfoQSTFT$ to nearby IF)} For all $m \in \{1, \dots, M\}$ and all $(t,\xi) \in \mathcal{Z}_m$ such that $|\QSTFT(t,\xi)| > \epstilde$, we have $|\IFinfoQSTFT(t,\xi) - \phi'_m(t)| \leq \epstilde$.  
\item \tbf{(Accuracy of reconstruction)} For each $m \in \{1, \dots, M\}$ there is a constant $C_m > 0$ where for any $t \in \R$, 
\begin{align*}
\left| \lim_{\beta \rar 0^+} \left(\int\limits_{ \{\xi \ : \ |\xi - \phi'_m(t)| < \epstilde \}} \dfrac{1}{\ \overline{h_{t,\xi}(0)} \ } \cdot \SSTQSTFTepstilde(t,\xi) \dd{\xi} \right) - f_m(t)\right| &\leq C_m \epstilde.
\end{align*}
\end{itemize}
\end{theorem}
\noindent We will give a proof of this theorem in future work.


\section{SST-QSTFT: discrete setting}
The discrete theory and implementation of SST-QSTFT differs depending on whether one uses a window that is band-limited, time-limited, or neither band-limited nor time-limited.  While the continuous theory prescribes the use of a band-limited window, it has nonetheless been proven\cite{behera2016theoretical} that the SST is still accurate for windows that are almost band-limited (i.e., having Fourier transform almost zero for all frequencies outside a certain passband). Moreover, for many signal processing applications, it is more common to use time-limited windows such as the well-known Hanning, Hamming, Blackman, and Kaiser windows. Indeed, it is generally not practical to use band-limited windows for real-time applications, since the implementation requires computing the DFT of the entire signal (which is not available).   Still, the usage of band-limited windows may enable sharper frequency resolution than is possible in the case of time-limited windows, which all exhibit some sort of frequency roll-off.  Moreover, certain windows may be neither band-limited nor time-limited, but are still of practical interest, such as the Gaussian window.  In practice, one truncates such windows and implements them in either the same way as a band-limited window or as a time-limited window.  

For the sake of brevity, we will focus on the most common case of time-limited windows, which enables extensions to the recently developed real-time SST theory\cite{su2016minimum, fourer2016recursive, chui2016realtime}. 

\subsection{Discrete quilted windows}\label{section-discrete-quilted-windows}
\begin{notation}
For $L \in \N$, we define $\Lset := \{0, \dots, L-1\}$. This notation will help to simplify the following definitions.
\end{notation}
\begin{defn}
We define a \emph{discrete quilted window family} for a given hop-size $H \in \N$, discrete-time signal of interest $f$, and maximal window length $L \in \N$, to be a tuple $\{\left(h_{n,k}, L_{n,k}\right)\}_{[n,k] \in \Nset \times \Lset}$, where for each $n \in \Nset$ and $k \in \Lset$,
\begin{itemize}
\item $h_{n,k}$ is a discrete sequence of length $L_{n,k} \in \N$ defining a discrete window function, 
\item the $L_{n,k}$ are chosen such that $L = \displaystyle\max_{[n,k] \in \Nset \times \Lset} L_{n,k}$, and
\item $N\in \N$ is large enough so that $f[(L-1) + (N-1)H] = 0$.
\end{itemize}
Hence, we associate to each time-frequency point $[n,k] \in \Nset \times \Lset$ a window function $h_{n,k}$ of length $L_{n,k}$.
\end{defn}\vskipp

\begin{remark}
It is well known that windows of different lengths have different time-frequency resolution. On the other hand, one often extends windows by zeros along their boundary to a length of $2^P$ for $P \in \N$ (a practice commonly known as \emph{zero-padding}) for the purpose of fast computation using the Fast Fourier Transform (FFT).  For our purposes, one should interpret the parameter $L_{n,k}$ that determines the length of each window in the discrete quilted window family as \emph{the window length after zero-padding}.
\end{remark}\vskipp

\begin{remark}
As in the continuous case, we do not generally choose an entirely different window for each time-frequency point, opting instead to divide the discrete time-frequency plane into tiles $T \subset \Nset \times \Lset$ , with a window $h^T$ associated to each $T$; that is, where $h_{n,k} = h^T$ for each $[n,k] \in T$. The constructions described in Sec.\ \ref{section-numerical-implementation-windows} proceed in such a manner.
\end{remark}

\subsection{Discrete quilted short-time Fourier transform (QSTFT)}
\begin{defn}
In analogy with the manner in which we defined the discrete STFT in Eqn.\ (\ref{discrete-STFT}), we may now define the \emph{discrete quilted short-time Fourier transform} of $f$  with \emph{hop size} $H \in \N$, for a discrete quilted window family $\{\left(h_{n,k}, L_{n,k}\right)\}_{[n,k] \in \Nset \times \Lset}$, via
\begin{equation}\label{discrete-QSTFT}
\QSTFTd[n,k] := \sum_{\ell=0}^{L_{n,k}-1} f[\ell + nH] \overline{h_{n,k}[\ell]} \e^{-\twopi\im k\ell/L_{n,k}},
\end{equation}
for each \emph{frame} $n \in \Nset$, and each \emph{frequency bin} $k \in \Lset$, where we define the discrete function $h$ by $h[\ell; n, k] := h_{n,k}[\ell]$ for each $n \in \Nset$ and $\ell, k \in \Lset$.
\end{defn}
 
\subsection{Discrete reassignment frequency formula}\label{section-alternative-reassignment-frequency}
In the continuous context, the reassignment frequency for SST-QSTFT requires the computation of the continuous-time derivative of the QSTFT $\QSTFT(t,\xi)$.  Since one cannot directly compute $\pt_t\QSTFT(t,\xi)$ in the discrete context, one option is to observe that $\pt_t\QSTFT(t,\xi)$ can be expressed as the QSTFT of $f$ with respect to the window function $(\twopi \im \xi h_{t,\xi}(\cdot-t) - \pt_t(h_{t,\xi})(\cdot-t))$. This is the method used by Thakur and Wu for the theory of the discrete SST-STFT\cite{thakur2011synchrosqueezing}. However, as observed in Sec.\ \ref{section-sst-stft-discrete-theory}, $\pt_t(h_{t,\xi})$ may be nonzero at the boundary of the support of $h_{t,\xi}$, which violates the periodicity requirement for taking the DFT and may cause aliasing.  Alternatively, we may wish to use $h_{n,k}$ that does not arise from any continuous-time function $h$, in which case such a quantity $\pt_t(h_{t,\xi})$ is unavailable.\vskipp

Instead, one may consider for instance a first-order finite difference approximation $D_{\Delta t}$ to the derivative $\pt_t$. However, using this approximation together with the continuous reassignment frequency formula will yield an inaccurate result in general, which we show in the following. First we put forth the following notation.

\begin{notation}
Suppose $s(\cdot)$ is a continuous-time signal defined on a set $\mathcal{R} \subset \R$, and we define the discretized version of $s$ by $s[n] := s(n\Delta t)$ for all $n \in \Z \cap \{n : \exists r \in \mathcal{R}, n = r/\Delta t\}$.  Then we define the \emph{backward-shifted discrete signal} $s^+$ by $s^+[n] := s((n+1)\Delta t)$.
\end{notation}\vskipp

\noindent The notation above becomes important because in general we may deal with a QSTFT having hopsize $H > 1$, and hence more coarsely sampled than the original signal. Now consider the approximation of the continuous reassignment frequency using the forward difference $(\QSTFTdp - \QSTFTd)/(\Delta t)$, for a complex-valued constant chirp of the form $f(t) = A \e^{\twopi\im c t}$, where $A, c \in \R^+$:
\begin{align*}
\dfrac{\pt_t \QSTFT(nH\Delta t, k/L_{n,k})}{\twopi\im\QSTFT(nH\Delta t, k/L_{n,k})} &\approx \dfrac{\QSTFTdp[n,k] - \QSTFTd[n,k]}{\twopi\im\Delta t\QSTFTd[n,k]}
\\&= \dfrac{ \sum_{\ell=0}^{(L_{n,k}-1)} A\left( \e^{\twopi\im c(\ell + nH + 1)\Delta t} - \e^{\twopi\im c(\ell + nH)\Delta t}\right) \overline{h_{n,k}[\ell]} \e^{-\twopi\im k\ell/L_{n,k}} }{\twopi\im\Delta t\sum_{\ell=0}^{(L_{n,k}-1)}  A\e^{\twopi\im c(\ell + nH)\Delta t} \overline{h_{n,k}[\ell]} \e^{-\twopi\im k\ell/L_{n,k}} }
\\&= \dfrac{\e^{\twopi\im c \Delta t} - 1}{\twopi\im\Delta t} = c \cdot \e^{\pi \im c \Delta t}\sinc(\pi c \Delta t).
\end{align*}
Hence, this forward difference does not isolate the IF $c$ in this simple case. Instead, one may consider using the alternative continuous reassignment frequency proposed by Oberlin et al.\cite{oberlin2014fourier}\ where one instead computes the time derivative of the \emph{phase spectrum} $\text{arg}(\QSTFT(t,\xi))$. Then, we see that the forward difference approximation exactly retrieves the IF $c$ for the complex-valued constant chirp $f$, provided that $c$ is below the Nyquist frequency $1/(2\Delta t)$:
\begin{align*}
\dfrac{1}{\twopi}\pt_t \text{arg}(\QSTFT(nH\Delta t, k/L_{n,k})) &\approx \dfrac{\text{arg}(\QSTFTdp[n,k]) - \text{arg}(\QSTFTd[n,k])}{\twopi\Delta t}
\\&= \dfrac{1}{\twopi\Delta t}\text{arg}\left(\dfrac{ \sum_{\ell=0}^{(L_{n,k}-1)} A\e^{\twopi\im c(\ell + nH + 1)\Delta t} \overline{h_{n,k}[\ell]} \e^{-\twopi\im k\ell/L_{n,k}} }{\sum_{\ell=0}^{(L_{n,k}-1)}  A\e^{\twopi\im c(\ell + nH)\Delta t} \overline{h_{n,k}[\ell]} \e^{-\twopi\im k\ell/L_{n,k}} }\right) = \dfrac{\text{arg}(\e^{\twopi\im c\Delta t})}{\twopi \Delta t} = c.
\end{align*}
Therefore, defining the \emph{discrete SST-QSTFT reassignment frequency}
\begin{align}\label{new-IF-estimate}
\IFinfoQSTFTd[n,k] &:= \dfrac{1}{2\pi \Delta t} \text{arg}\left( \dfrac{\QSTFTdp[n,k]}{\QSTFTd[n,k]} \right)
\end{align}
for each frame $n \in \Nset$ and frequency bin $k \in \Lset$ where $\QSTFTd[n,k] \neq 0$, we see that $\IFinfoQSTFTd[n,k]$ is \emph{exact} for the complex-valued constant chirp of frequency $c \in (0, \frac{1}{2\Delta t})$ considered.  Moreover, the derivation of the continuous formula for the reassignment frequency is motivated by the fact that it is exact for a constant chirp \cite{daubechies2011synchrosqueezed, daubechies1996nonlinear}. Hence, we use the discrete reassignment formula given by (\ref{new-IF-estimate}) in the following definition of the discrete version of SST-QSTFT.
\subsection{Discrete QSTFT-based SST (SST-QSTFT)}\label{section-discrete-sst-qstft}
\begin{defn}
We define the \emph{discrete QSTFT-based synchrosqueezing transform} (discrete SST-QSTFT)
of a discrete-time signal $f$, with respect to the function $h$ defining the discrete quilted window sequences $\{h_{n,\ell}\}_{[n,\ell] \in \Nset \times \Lset}$, and with tolerance $\gamma\geq 0$, as follows:
\begin{displaymath}\label{equation-sst-qstft-discrete-original}
\SSTQSTFTd[n,k] := \sum_{\ell \in \AsetQSTFTd} \QSTFTd[n,\ell] \charone_{\BsetQSTFTd}[\ell],
\end{displaymath}
for each frame $n \in \Nset$ and for each \emph{reassignment frequency bin} $k \in \Kset$ with $K \geq L$, where $\AsetQSTFTd := \{\ell \in \Lset : |\QSTFTd[n,\ell]| > \gamma\}$, $\BsetQSTFTd := \{ \ell \in \Lset \ : \ -\frac{1}{2} \leq K\Delta t\IFinfoQSTFTd[n,\ell] - k < \frac{1}{2}\}$ is the set of QSTFT frequency bins $\ell$ at the frame $n$ where the corresponding reassignment frequency is closer to $k$ than any other reassignment frequency bin, $\charone_X$ is the characteristic function on the set $X$ (i.e., $\charone_X(x) = 1$ if $x \in X$ and $\charone_X(x) = 0$ otherwise), and $\IFinfoQSTFTd$ as defined in (\ref{new-IF-estimate}) is the \emph{discrete QSTFT-based reassignment frequency.}
\end{defn}\vskipp

\begin{remark}\label{remark-no-b-approximate-summation}
Here, we do not use any function $b$ to do an approximate summation, as in the case of the continuous SST-QSTFT. The use of $b$ in the continuous case is for the sake of the proof of reconstruction, and in the following we do not provide a direct reconstruction formula from the discrete SST-QSTFT coefficients, passing back instead to the original QSTFT coefficients. See Sec.\ \ref{section-discrete-time-limited-reconstruction} for more details.
\end{remark}\vskipp

\begin{remark}
Note that the SST-QSTFT may be computed over more frequency bins $K$ than the QSTFT has ($L$), enabling the possibility of an even more precise estimation of IF. Indeed, since the discrete reassignment frequency formula is exact in the case of a constant chirp \emph{regardless of whether the IF coincides with a frequency bin,} it follows that the use of $K > L$ may well yield a more precise concentration around the true IF.
\end{remark}\vskipp

\begin{defn}
Alternatively, we may also calculate the \emph{discrete QSTFT-based magnitude synchrosqueezing transform} (discrete magnitude SST-QSTFT), where we reassign the \emph{magnitude-squared} of the QSTFT coefficients, by
\begin{displaymath}
\SSTQSTFTdMAG[n,k] := \sum_{\ell \in \AsetQSTFTd} \left|\QSTFTd[n,\ell]\right|^2 \charone_{\BsetQSTFTd}[\ell].
\end{displaymath}
\end{defn}\vskipp

\begin{remark}\label{remark-sstqstftdmag}
Calculating $\SSTQSTFTdMAG$ enables us to preserve the total energy of the QSTFT coefficients, and hence yields higher energy peaks along IF curves than the original QSTFT.  Since the original QSTFT coefficients are summed together for $\SSTQSTFTd$, their total energy might decrease after summation.  Both SST-QSTFT representations will be more concentrated than QSTFT, due to the squeezing procedure summing together high-energy coefficients closer to the IF curves, but $\SSTQSTFTdMAG$ will generally contain more energy along the IF curves than $\SSTQSTFTd$.
\end{remark}\vskipp 

\begin{remark}
One may set $\gamma > 0$ to ensure an accurate reassignment frequency for all reassigned coefficients.  In fact, this is necessary in the continuous case to guarantee reconstruction accuracy, as per Theorem \ref{theorem-continuous}.  However, as stated in Remark \ref{remark-no-b-approximate-summation}, we do not perform direct reconstruction from the discrete SST-QSTFT coefficients.  Moreover, we use $\SSTQSTFTdMAG$ in all our numerical experiments (for reasons given in Sec.\ \ref{section-discrete-time-limited-reconstruction}), and the issue of inaccuracy for reassignment frequency only concerns low-energy coefficients that do not contribute much to $\SSTQSTFTdMAG$. Hence, we always set $\gamma = 0$ in our experiments.
\end{remark}

\subsection{Adaptive discrete quilted window function families}
For the case of discrete quilted window families, we make the following definition of an adaptive family:
\begin{defn}
Suppose that $f \in \Bepsd$. We say that the quilted window family $\{\left(h_{n,k}, L_{n,k}\right)\}_{[n,k]\in\Nset\times \Lset}$ is of class $W^Q_{d,\eps,f}$ if the following conditions hold: 
\begin{itemize}
\item \tbf{Fourier-side decay (almost band-limitation):} There exists a constant $C_h \in \R^+$ such that
\begin{align}\label{definition-Ch}
\displaystyle\sup_{\substack{u\in\R \setminus [-d/2,d/2] \\ [n,k] \in \Nset \times \Lset}}\left|\widehat{h_{n,k}}(u)\right| = C_h\eps,
\end{align} 
so that the semi-discrete Fourier transform of $h_{n,k}$ is bounded outside the frequency band $[-d/2, d/2]$. 
\item \tbf{Window choice remains constant in the frequency band around an IF value:} For all $n \in {\Nset}$ and $m \in \{1, \dots, M\}$, there exists a single discrete window function $g_{n,m}$ such that $h_{n,k} \equiv g_{n,m}$ for all $k \in \{k : |k/(L_{n,k}\Delta t) - \phi'_m[nH]| < d/2\}$.
\end{itemize}
Moreover, we call $W^Q_{d,\eps,f}$ the \emph{class of $f$-adaptive discrete quilted window families}. 
\end{defn}

\subsection{Theorem}
Now, we state the following theorem for the discrete SST-QSTFT:
\begin{theorem}\label{theorem-discrete-time-limited}
Let $\eps > 0$, $\nu \in (0,1/2)$, $\epstilde := \eps^{\nu}$, $d > 0$. Suppose that $\displaystyle f = \sum_{m=1}^M f_m \in \Bepsd$, with $\phi'_m[nH] \leq 1/(2\Delta t)$ for all $m \in \{1, \dots, M\}$ and $n \in \Nset$. Assume that $\{\left(h_{n,k}, L_{n,k}\right)\}_{[n,k] \in \Nset \times \Lset}$ is of the class $W^Q_{d,\eps,f}$. Then, if $\eps$ is sufficiently small we have:
\begin{itemize}
\item \tbf{(Concentration of QSTFT around IF curves)} $|\QSTFTd[n,k] | > \epstilde$ only when there is an $m \in \{1, \dots, M\}$ such that $[n,k] \in Z_m := \{ [n,k] \in \Nset \times \Lset:\ |\phi'_m[nH] - k/(L_{n,k}\Delta t)| < d/2 \}$.
\item \tbf{(Closeness of reassignment frequency $\IFinfoQSTFTd$ to nearby IF)} For all $m \in \{1, \dots, M\}$ and all $[n,k] \in Z_m$ such that $|\QSTFTd[n,k]| > \epstilde$, we have $|\IFinfoQSTFTd[n,k] - \phi'_m[nH]| \leq C_{\epstilde}$, where $\lim_{\epstilde \rightarrow 0^+} C_{\epstilde} = 0$.  
\end{itemize}
\end{theorem}
\noindent The proof of this theorem is left for our future work, due to the page limitation.

\subsection{Discrete reconstruction}\label{section-discrete-time-limited-reconstruction}
We do not provide a reconstruction theorem for the discrete case, because the technique of summing over SST-QSTFT reassignment frequency bins for each frame is not accurate in the case of time-limited windows used with hop size $H > 1$.  Instead, we refer to the technique of Holighaus et al.\cite{holighaus2016reassignment}, who suggested to store the \emph{inverse reassignment map} defined by 
\begin{align*}
\IFinfoQSTFTdInverse[n,k] &:= \{ \ell \in \{0, \dots, L-1\} \ : \IFinfoQSTFTd[n,\ell] = k \}
\end{align*}
for $n \in \Nset$ and $k \in \Kset$.  $\IFinfoQSTFTdInverse[n,k]$ is the set of all frequencies $\ell \in \Lset$ whose reassignment frequency is $k$. Then, rather than doing reconstruction along the narrow SST IF ridges, the reconstruction can be done over the thicker QSTFT ridges via frame synthesis, using an \emph{overlap-add formula} \cite{SASPWEB2011} and the construction of dual windows using canonical tight frames \cite[Sec.\ 3.3, Theorem 1]{balazs2011theory}. Note that for this type of reconstruction, the condition of maintaining the same quilted window over all frequencies in a given IF band (for a fixed frame) is essential.

We note that the usual technique for reconstruction using SST is to: (i) analyze the magnitude-squared of the SST representation for peaks that form continuous ridges; (ii) use an algorithm to extract IF curves following these ridges\cite{auger2013time, iatsenko2013extraction, meignen2017fully}; and then (iii) employ a discrete version of a reconstruction formula such as the one in Theorem \ref{theorem-continuous} to reconstruct from the original SST coefficients\cite{daubechies2011synchrosqueezed, oberlin2014fourier}.  By contrast, the technique above enables us to replace step (iii) with reconstruction from QSTFT coefficients, and we only need the SST for extracting the IF curves.  Hence, rather than using $\SSTQSTFTd$, it makes more sense to use the generally more strongly concentrated $\SSTQSTFTdMAG$ for analysis (see Remark \ref{remark-sstqstftdmag}). \emph{In the following, we will refer to both $\SSTQSTFTd$ and $\SSTQSTFTdMAG$ as SST-QSTFT, and use only $\SSTQSTFTdMAG$ in our numerical experiments.} 

\section{Numerical implementation}\label{section-numerical-implementation}
In the following section, we describe the numerical implementation of SST-QSTFT.  We have developed a Python suite for the adaptive time-frequency transforms described in this work, including SST-QSTFT, available to researchers upon request.\footnote{Please email Alex Berrian at \url{aberrian@math.ucdavis.edu} for a copy of the Python suite.}

\subsection{Automatic adaptive window selection}\label{section-numerical-implementation-windows}
The question remains how to adapt the window selection to the signal. In 2007, Jaillet \& Torr\'{e}sani introduced the idea of the \emph{time-frequency jigsaw puzzle}, where different windows could be adaptively associated to different regions (representing ``jigsaw pieces'') in the time-frequency plane, based on the time-frequency content \cite{jaillet2007time}.  As a starting point, we may incorporate only a single iteration of their algorithm in order to choose an optimal window for each time-frequency region, from a given window collection.  Our first algorithm for automatic adaptive window selection proceeds as follows:

\begin{algorithm}\label{algorithm-1}
\begin{itemize}
\item[\tbf{1.}] Choose windows $h^w$, $w = 1, \dots, W$ for $W \in \N$, and corresponding sampling lattices, the latter determined by the hop size $H$ (which must be the same for all windows, by Sec.\ \ref{section-discrete-quilted-windows}) and FFT size $L^w$.
\item[\tbf{2.}] For each $h^w$, compute the discrete STFT $V_{h^w} f$ with hop size $H$ and FFT size $L^w$ (chosen in Step 1).  
\item[\tbf{3.}] Fix parameters $A, B \in \N$ to decompose the time-frequency plane $\R^2$ into ``supertiles'' $\Box_{r,s}$, where \\$\Box_{r,s} := [r\tilde{A}, (r + 1)\tilde{A}) \times [s\tilde{B}, (s+1)\tilde{B})$, $\tilde{A} := A H\Delta t$, $\tilde{B} := B/(F_\text{min} \Delta t)$, and $F_\text{min} := \text{min}_w \{F^w\}$.  We note that $A$ (resp.\ $B$) represents the number of points on the coarsest grid along the time (resp.\ frequency) axis contained in each supertile.
\item[\tbf{4.}] For each $\Box_{r,s}$ and $h^w$, calculate the total energy $E^w_{r,s} := \sum_{(\ell,k) \in \Box_{r,s}}|\QSTFTdhw[\ell,k]|^2$.
\item[\tbf{5.}] Fix $\alpha \in (0,1)$. Then for each $\Box_{r,s}$ and $h^w$, calculate the \emph{sampled R\'{e}nyi entropy}\cite{liuni2013thesis, jaillet2007time} $\mathcal{R}^{\alpha,w}_{r,s}$ of the STFT coefficients whose corresponding lattice points are within $\Box_{r,s}$:
\begin{align*}
\mathcal{R}^{\alpha,w}_{r,s} &:= \dfrac{1}{1-\alpha}\log \left[\displaystyle\sum_{(\ell,k) \in \Box_{r,s}} \left(\left(\frac{H}{L^w}\right)^{1-\alpha}\frac{|\QSTFTdhw[\ell,k]|^2}{E^w_{r,s}} \right)^\alpha\right] 
\end{align*}
\item[\tbf{6.}] To each supertile $\Box_{r,s}$, associate the window $h^{w^*}$ with the smallest sampled R\'{e}nyi entropy $\mathcal{R}^{\alpha,w^*}_{r,s}$.  Thus we adapt to signal content in $\Box_{r,s}$.
\end{itemize}
\end{algorithm}
\begin{remark}
The R\'{e}nyi entropy can be seen as a generalization of the $\ell^2$-normalized $\ell^1$-norm, and can hence be seen as a time-frequency sparsity measure. Indeed, with $\alpha = 0.5$, one recovers the $\ell^2$-normalized $\ell^1$-norm.  Furthermore, minimization of the $\ell^2$-normalized $\ell^p$-norm with $p = 2\alpha$ is equivalent to minimizing the R\'{e}nyi entropy\cite{ricaud2014survey}. 
\end{remark}\vskipp

\begin{figure}
\begin{center}
\includegraphics[width=18cm, height=6cm]{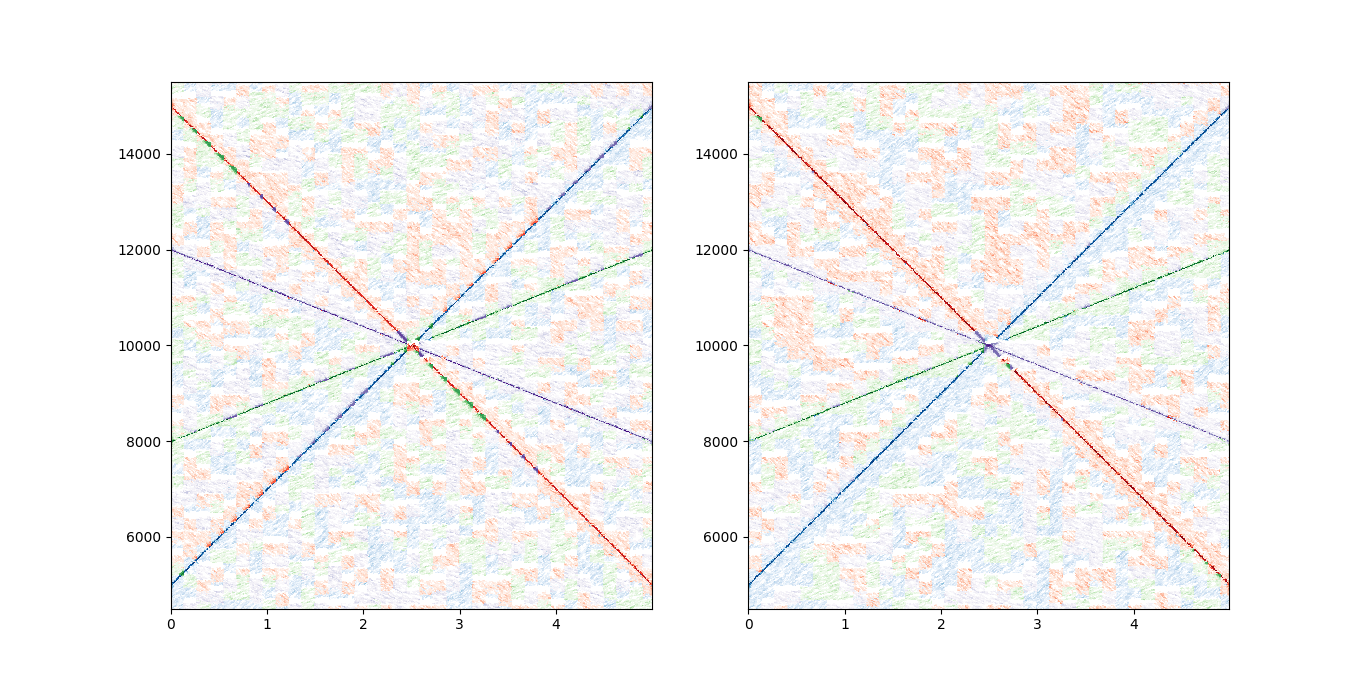}
\end{center}
\caption{\label{figure-qstft-sst-qstft-nonperturbed} \tbf{Left:} SST-QSTFT of the crossing chirps signal $f$ using Algorithm \ref{algorithm-1}.  \tbf{Right:} SST-QSTFT of the same signal using Algorithm \ref{algorithm-2}.  In both plots, colors correspond to the window chosen in each supertile for the calculation of the QSTFT. Blue: $h^1$. Green: $h^2$. Purple: $h^3$. Red: $h^4$. Frequency range is restricted from 4500 Hz to 15500 Hz. The perturbed supertiles algorithm leads to fewer instances of QSTFT coefficients near IFs from windows of very dissimilar chirp rate.}
\end{figure}
The left plot in Figure \ref{figure-qstft-sst-qstft-nonperturbed} demonstrates the application of Algorithm \ref{algorithm-1} to compute the SST-QSTFT of a noise-contaminated signal $f(t) = \sum_{m=1}^4 f_m(t) + f_\text{noise}(t)$ containing crossing linear chirp components of the form $f_m(t) = \cos(\twopi \phi_m(t))$ with $\phi'_m(t) = c_m + \sigma_m t$, where $(c_m)_{m=1}^4 = (5000, 8000, 12000, 15000)$ and $(\sigma_m)_{m=1}^4 = (2000, 800, -800, -2000)$, and with zero-mean white noise component $f_\text{noise}(t)$ such that $f$ has SNR equal to $4.0$ dB.  Here, SNR for a generic signal $f(t) = y(t) + f_\text{noise}(t)$, where $y$ is noiseless, is defined by \cite{daubechies2011synchrosqueezed}
\begin{align*}
\text{SNR}(f) &:= 10\log_{10}\left( \dfrac{\text{var}\left(y\right)}{\text{var}\left(f_\text{noise}\right)}\right).
\end{align*}
To compute the SST-QSTFT, we use the chirped window collection $\{h^w\}$ where $h^w(t) := h^0(t)\e^{\twopi \im \sigma^w t^2/2}$.  Here $h^0$ is a Hanning window of width $4000$ samples, zero-padded to $L = 2^{12} = 4096$ samples, and $\{\sigma^w\}_{w=1}^4 = \{ 1900, 900, -900, -1900\} \text{Hz}^2$.  We set $H =250$, $K = 4L = 16384$, $\gamma = 0$, $\alpha = 0.5$, and $A = B = 24$.\vskipp

\begin{remark}
We use this particular set of chirp rates $\{\sigma^w\}$ only because it demonstrates that a reasonable result may be obtained using windows with chirp rates that are close, but not quite equal, to the actual local chirp rates $\{\sigma_m\}$.  Often, $\{\sigma^w\}$ may be chosen by visual inspection of a single STFT of the signal.  However, in general, the local chirp rates are unknown, and one should have some systematic method of deciding the collection $\{\sigma^w\}$.  This topic is out of the scope of this paper, and we leave it for our future work.
\end{remark}\vskipp

In Figure \ref{figure-qstft-sst-qstft-nonperturbed}, we see that supertiles in the neighborhood of the chirps may be assigned a window with chirp parameter that does not correspond to the slope of the nearby chirp component.  This may happen in the case that the supertile contains very little signal content, in which case sharply-concentrated content may be mistaken for noise.  There is also the possibility that an IF curve does not pass through a supertile at all, but is analyzed by a window that yields a very blurry representation of the IF, which ends up ``leaking'' into the supertile.  In both of these cases, Algorithm \ref{algorithm-1} may fail to pick the window that yields the sparsest signal representation in the supertile.

\begin{figure}
  \begin{minipage}[c]{0.65\textwidth}
    \includegraphics[width=\textwidth, height=0.9\textwidth]{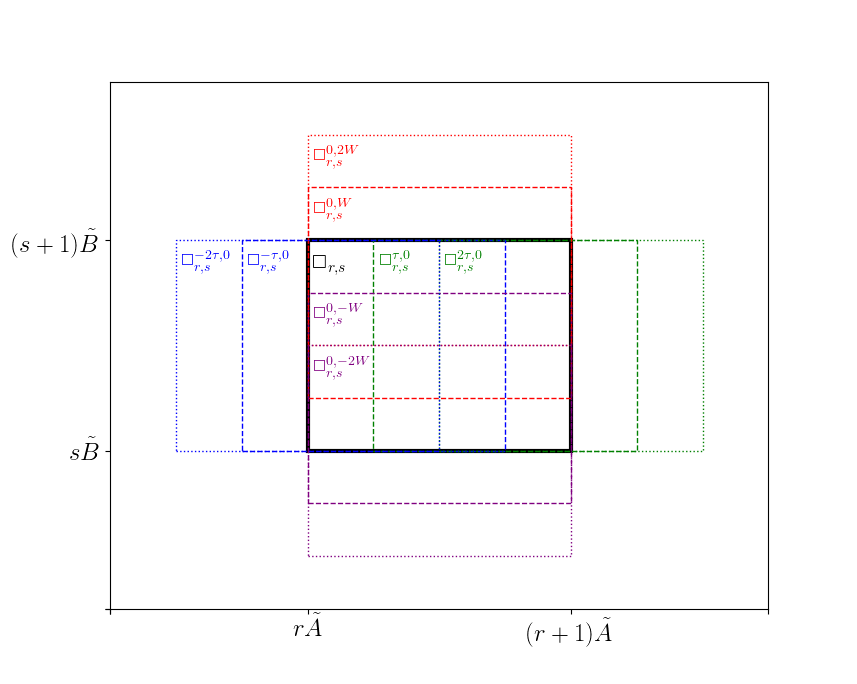}
  \end{minipage}\hfill
  \begin{minipage}[c]{0.35\textwidth}
    \caption{\label{figure-supertile-perturbations} All perturbations of the supertile $\Box_{r,s}$, for Algorithm \ref{algorithm-2} parameters given as follows: $P_T = P_Y = \{-1,1\}$ (backward and forward directional perturbation in both time and frequency), $\Tshift = \Yshift = 2$ (two perturbations per direction in both time and frequency), $\tau = \Tstep \cdot ( H\Delta t)$, and $W = \Ystep/(F_\text{min} \Delta t)$ for arbitrary $\Delta t$, $H, F_\text{min}, \Tstep$, and $\Ystep$.  The label $\Box^{\tilde{t}, \tilde{y}}_{r,s}$ of each perturbed supertile is shown in the upper-left corner of that supertile.}
 \end{minipage}
\end{figure}

We note that Algorithm \ref{algorithm-1} is essentially equivalent to a simplified version of the algorithm of Sheu et al.\cite{sheu2015entropy}, but with the R\'{e}nyi entropy written out explicitly for the discrete setting.  Sheu et al.\ allow for overlapping supertiles, and update the optimal entropy choice as they move forward in time and frequency.  To improve upon the deficiencies of Algorithm \ref{algorithm-1}, we go with a more general approach based on \emph{perturbed supertiles}.  In Figure \ref{figure-supertile-perturbations}, we plot a single supertile $\Box_{r,s}$ together with eight perturbations of $\Box_{r,s}$ along the time and frequency directions.  The basic idea of Algorithm \ref{algorithm-2} is to calculate the R\'{e}nyi entropy in the original supertile $\Box_{r,s}$ as well as in several perturbations $\Box_{r,s}^{\tilde{t}, \tilde{y}}$ of $\Box_{r,s}$, and then to compute the average $\widetilde{\mathcal{R}}^{\alpha, w}_{r,s}$ of all the calculated entropies, for all windows (indexed by $w$).  Then, the window that minimizes this averaged entropy $\widetilde{\mathcal{R}}^{\alpha, w}_{r,s}$ is the optimal window for the supertile $\Box_{r,s}$.  We give a rigorous description of our algorithm as follows:\vskipp

\begin{algorithm}\label{algorithm-2}
\begin{itemize}
\item[\tbf{1--3.}] Do Steps 1-3 of Algorithm \ref{algorithm-1}.
\item[\tbf{4.}] Fix permutation step sizes $\Tstep$ and $\Ystep$, the total amounts of permutations $\Tshift$ and $\Yshift$, and the sets $P_T$ and $P_Y$ of permutation directions satisfying $P_T, P_Y \subseteq \{-1, 1\}$, for the time ($T$) and frequency ($Y$) axes respectively.  We explain these parameters as follows:
\begin{itemize}
\item $\Tstep$ (resp.\ $\Ystep$) is given in units of points on the coarsest grid in the time (resp.\ frequency) direction.  For example, suppose one sets $\Ystep = 3$.  Then, since the distance between two points on the coarsest frequency grid is $1/(F_\text{min} \Delta t)$ Hz, one permutes the supertile with a step size of $W := 3/(F_\text{min}\Delta t)$ Hz.  If one sets $\Yshift = 2$, then the supertile is permuted two separate times in the frequency direction, each time by $W$ Hz (see Figure \ref{figure-supertile-perturbations}).
\item The possible forms of $P_T$ and $P_Y$ are $\{-1\}$, $\{1\}$, and $\{-1, 1\}$.  The first option means that the permutations of the supertile grid will be in the backward direction. The second option refers to the forward direction, and the last option refers to permutation in both directions.
\end{itemize}
\item[\tbf{5.}] Define $\widetilde{T} := \left\{ p_T \cdot t \cdot \Tstep \cdot ( H\Delta t): \ p_T \in P_T; \ t = 1, \dots, \Tshift \right\}$, and \\$\widetilde{Y} := \left\{ p_Y \cdot y \cdot \Ystep /(F_\text{min} \Delta t) : \ p_Y \in P_Y; \ y = 1, \dots, \Yshift \right\}$. Then: 
\begin{itemize}
\item[\tbf{5a.}] For each $\tilde{t} \in \widetilde{T}$, define $\Box^{\tilde{t}, 0}_{r,s} := [r\tilde{A} + \tilde{t}, (r + 1)\tilde{A} + \tilde{t}) \times [s\tilde{B}, (s+1)\tilde{B})$.
\item[\tbf{5b.}] For each $\tilde{y} \in \widetilde{Y}$, define $\Box^{0, \tilde{y}}_{r,s} := [r\tilde{A} , (r + 1)\tilde{A}) \times [s\tilde{B} + \tilde{y}, (s+1)\tilde{B} + \tilde{y})$.
\end{itemize}
\item[\tbf{6.}] Now, fix the R\'{e}nyi entropy parameter $\alpha \in (0,1]$ and do the following for each pair $(r,s)$: 
\begin{itemize}
\item[\tbf{6a.}] For each $\Box \in \mathfrak{B}_{r,s} := \left\{ \Box_{r,s}, \Box^{\tilde{t}, 0}_{r,s}, \Box^{0, \tilde{y}}_{r,s}: \ \tilde{t} \in \widetilde{T}; \ \tilde{y} \in \widetilde{Y} \right\}$ and $h^w$, calculate the total energy\\ $E^{w,\Box}_{r,s} := \displaystyle\sum_{(\ell,k) \in \Box}|\QSTFTdhw[\ell,k]|^2$.
\item[\tbf{6b.}] For each $\Box \in \mathfrak{B}_{r,s}$ and $h^w$, calculate the sampled R\'{e}nyi entropy $\mathcal{R}^{\alpha,w,\Box}_{r,s}$ of the STFT coefficients whose corresponding lattice points are within $\Box$:
\begin{align*}
\mathcal{R}^{\alpha,w,\Box}_{r,s} &:= \dfrac{1}{1-\alpha}\log \left[\displaystyle\sum_{(\ell,k) \in \Box} \left(\left(\frac{H}{L^w}\right)^{1-\alpha}\frac{|\QSTFTdhw[\ell,k]|^2}{E^{w,\Box}_{r,s}} \right)^\alpha\right]
\end{align*}
\item[\tbf{6c.}] For each $h^w$, calculate the averaged entropy measure 
\begin{align*}
\widetilde{\mathcal{R}}^{\alpha, w}_{r,s} &:= \frac{1}{1 + \left|\widetilde{T}\right| + \left|\widetilde{Y}\right|} \sum_{\Box \in \mathfrak{B}_{r,s}} \mathcal{R}^{\alpha,w,\Box}_{r,s}.
\end{align*}
\end{itemize} 
\item[\tbf{7.}] To each supertile $\Box_{r,s}$, associate the window $h^{w^*}$ with the smallest averaged entropy $\widetilde{\mathcal{R}}^{\alpha,w^*}_{r,s}$.  Thus we adapt to signal content in $\Box_{r,s}$.
\end{itemize}
\end{algorithm}
\begin{remark}
Algorithm \ref{algorithm-2} is significantly different from the overlapping supertiles algorithm of Sheu et al., since their algorithm only proceeds forward in time and frequency, and does not incorporate an averaging procedure.
\end{remark}\vskipp

\begin{remark}
Algorithm \ref{algorithm-2} is similar to a method used by Coifman and Donoho in the context of wavelet-based signal denoising, called ``cycle-spinning,'' that involved averaging a certain quantity over grid perturbations\cite{coifman1995translation}.
\end{remark}\vskipp

The intended effect of averaging the entropy over perturbations of the given supertile is to avoid the errors that result from supertiles containing minimal signal content and from blurred signal representation leaking into surrounding supertiles.  In the right plot of Figure \ref{figure-qstft-sst-qstft-nonperturbed}, we used Algorithm \ref{algorithm-2} with $\alpha=0.5$, $A = B = 24$, $\Tstep = \Ystep = \Tshift = \Yshift = 4$, and $P_T = P_Y = \{-1, 1\}$. We see that Algorithm \ref{algorithm-2} improves the performance of SST-QSTFT, leading to far fewer blurry patches along the IF curves.  This sharper representation results from the more frequent selection of windows with chirp parameter closest to the nearest IF component's chirp rate.\vskipp

\begin{remark}
We note that the selection of the parameters $\Tstep,$ $\Ystep,$ $\Tshift,$ and $\Yshift$ for Figure \ref{figure-qstft-sst-qstft-nonperturbed} (and later on for Figures \ref{figure-four-transforms} and \ref{figure-gibbons}) is somewhat arbitrary.  In general, some tuning of these parameters is necessary to achieve a desired result.  We do not address optimal strategies for tuning these parameters in this paper, leaving this topic for our future work.  However, we generally recommend using $P_T = P_Y = \{-1, 1\}$, so that the entropy averaging calculation is not biased either forward or backward in time or frequency.
\end{remark}

\begin{table}
\begin{center}
  \begin{tabular}{|c|c|c|c||c|c|c||c|c|c||c|c|c|}
    \cline{2-13}
     \multicolumn{1}{c|}{} & \multicolumn{3}{|c||}{1st-order SST-STFT} & \multicolumn{3}{|c||}{2nd-order SST-STFT} & \multicolumn{3}{|c||}{RM}& \multicolumn{3}{c|}{SST-QSTFT} \\
     \cline{2-13}
      \multicolumn{1}{c|}{} & 6 dB & 0 dB & -6 dB & 6 dB & 0 dB & -6 dB & 6 dB & 0 dB & -6 dB & 6 dB & 0 dB & -6 dB \\ \hline
    $\phi'_1$ & 3.3\%  & 2.6\% & 1.5\% & 9.7\%  & 5.4\%  & 2.1\% & 15.4\% & 10.8\% & 4.6\% & 16.9\% & 13.0\% & 6.3\% \\\hline
    $\phi'_2$ & 10.5\% & 7.9\% & 4.0\% & 19.8\% & 13.7\% & 6.0\% & 21.8\% & 16.1\% & 7.6\% & 20.6\% & 15.7\% & 8.1\% \\\hline
    $\phi'_3$ & 10.6\% & 8.4\% & 4.4\% & 20.2\% & 14.2\% & 6.5\% & 22.2\% & 16.7\% & 8.2\% & 21.0\% & 16.2\% & 8.9\% \\\hline
    $\phi'_4$ & 3.3\%  & 2.6\% & 1.4\% & 9.9\%  & 5.3\%  & 2.0\% & 15.3\% & 10.6\% & 4.5\% & 17.1\% & 12.8\% & 6.1\% \\\hline
  \end{tabular}\vspace*{3mm}\\
  \caption{Comparison of percent total energy along IF ridges relative to the total energy between 5 and 15 kHz over all times, for time-frequency representations 1st-order SST, 2nd-order SST, RM, and SST-QSTFT, and SNRs of 6 dB, 0 dB, and -6 dB.  Specifically, we compute $E_C(R^q_m(C,f))/E_C(\Z_N \times \mathfrak{W})$ for each time-frequency representation $C$, where $E_C(R)$ denotes $\sum_{[n,k] \in R} C[n,k]$ over the given region $R$, $f$ is the crossing chirps signal, and $\mathfrak{W}$ denotes bins between 5 and 15 kHz.}
  \label{table-second-order-sst-vs-sst-qstft}
\end{center}
\end{table}

\subsection{Numerical results}

\subsubsection{Comparison with SST-STFT, second-order SST, and reassignment method}
One way to measure the sharpness of a time-frequency representation for a given signal is to calculate the amount of energy along the \emph{ridges} of the IF components.  Given a signal $f = \sum_{m=1}^M f_m \in \Bepsd$ and a discrete time-frequency representation $C(f)$ (for instance, STFT or SST-STFT) operating on $f$ with $N$ frames, $K$ frequency bins, and hop size $H$, we define the $m^\text{th}$ ridge ($m = 1, \dots, M$) of bandwidth $q$ bins to be $R^q_m(C,f) := \{ C(f)[n,k] : \ -\frac{1}{2} - q \leq K\Delta t\phi'_m[nH] - k < \frac{1}{2} + q\}  \}$; i.e., all the coefficients $C(f)[n,k]$ within $q$ frequency bins from the closest frequency bin to $\phi'_m[nH]$.  

Using the crossing chirps signal, we compare the ridge concentration performance of SST-QSTFT with SST-STFT, \emph{second-order SST-STFT}\cite{oberlin2015second}, and the \emph{reassignment method} (RM)\cite{kodera1976new, auger1995improving}.  The central idea of second-order SST-STFT is to use an improved, ``second-order accurate'' reassignment frequency formula that exactly retrieves the IF in the case of a noiseless single-component signal with IF of the form $\phi'(t) = \sigma t + c$ for real-valued constants $\sigma$ and $c$.  Second-order SST-STFT requires the side computation of a \emph{reassignment time} quantity, but still only reassigns the frequency locations of the STFT coefficients.  By contrast, RM uses reassignment times together with reassignment frequencies to reassign \emph{both} the time and frequency locations of the STFT coefficients.  

In order to compute second-order SST-STFT and RM, we introduce the \emph{discrete reassignment time formula}
\begin{align}\label{definition-RT}
\RTSTFTd[n,k] &:= nH\Delta t - \dfrac{1}{2\pi} \text{arg}\left( \dfrac{\STFTdfp[n,k]}{\STFTd[n,k]} \right),
\end{align}
where $\STFTdfp[n,k] := \sum_{\ell=0}^{L-1} f[\ell + nH] \overline{g[\ell]} \e^{-\twopi \im \ell(k/L + \Delta t)}$ is the STFT of $f$ with normalized frequency shifted ahead by $\Delta t$.  It is easy to show that for a single spike $f(t) = A \cdot \charone_{\{t = \ell_0\Delta t\}}(t)$ with constant amplitude $A > 0$ located at $t = \ell_0 \Delta t$, $\RTSTFTd[n,k] \equiv \ell_0 \Delta t$ whenever $nH - \ell_0 \in \Lset$.  This formula is a discretization of the continuous reassignment time formula given by $\RTSTFT(t,\xi) := t - \frac{1}{\twopi}\pt_\xi \text{arg}\left(\STFT(t,\xi)\right)$\cite{oberlin2015second}.  For second-order SST-STFT specifically, we introduce the \emph{discrete second-order reassignment frequency formula} given by
\begin{align}\label{definition-second-order-discrete-RF}
\IFinfoSTFTdsecond[n,k] &:= \IFinfoSTFTd[n,k] + \dfrac{D_t\IFinfoSTFTd[n,k]}{D_t\RTSTFTd[n,k]}\left( (nH + \lfloor0.5L\rfloor)\Delta t - \RTSTFTd[n,k]\right),
\end{align}
where $D_t\IFinfoSTFTd[n,k] := \frac{1}{2\Delta t} \left( \IFinfoSTFTdfplus[n,k] - \IFinfoSTFTdfminus[n,k]\right)$, $D_t\RTSTFTd[n,k] := \frac{1}{2\Delta t} \left( \RTSTFTdfplus[n,k] - \RTSTFTdfminus[n,k]\right)$, and $f^-[\ell] := f((\ell-1)(\Delta t))$ is the \emph{forward-shifted discrete signal}.  This is a discretization of the continuous formulation of second-order reassignment frequency given in Eqn.\ (32) in the seminal paper by Oberlin et al.\cite{oberlin2015second}\ on second-order SST.  We use this discretization instead of the alternative formulation proposed by those authors in Eqns.\ (13) and (31) of the aforementioned paper, because it allows us to avoid the issues mentioned in Sec.\ \ref{section-alternative-reassignment-frequency}.

\begin{figure}
\begin{center}
\includegraphics[width=18cm, height=6cm]{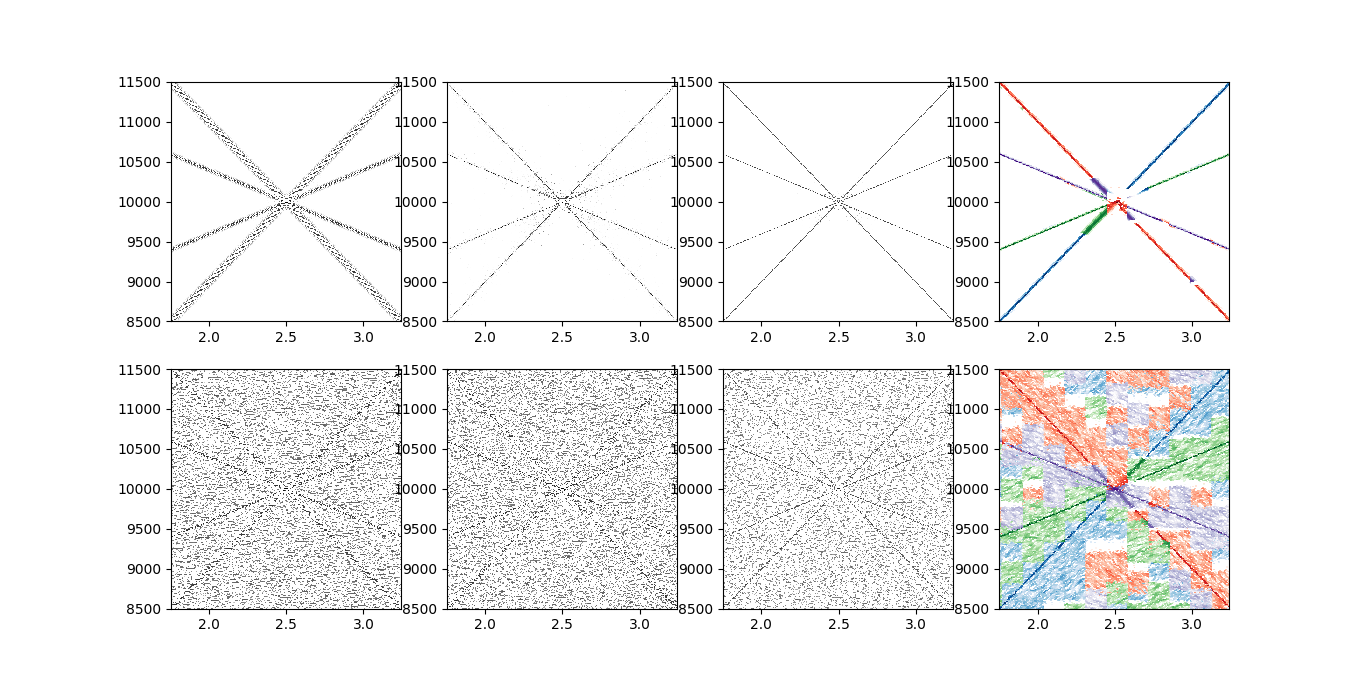}
\end{center}
\caption{\label{figure-four-transforms} Time-frequency representations of noiseless \tbf{(top row)} and -6 dB SNR \tbf{(bottom row)} crossing chirps signal.  \tbf{From left to right:} 1st-order SST, 2nd-order SST, RM, and SST-QSTFT.  For visualization purposes, we focus on the time-frequency region between $1.5$ and $2.5$ seconds and $8.5$ and $11.5$ kHz.  In SST-QSTFT plots, colors correspond to the window chosen in each supertile for the calculation of the QSTFT. Blue: $h^1$. Green: $h^2$. Purple: $h^3$. Red: $h^4$.  The SST-QSTFT enables the IFs to be more easily visualized under noise, and concentrates the most signal energy along the IF ridges, as per Table \ref{table-second-order-sst-vs-sst-qstft}. \emph{Note:} The calculation in Table \ref{table-second-order-sst-vs-sst-qstft} is done over a larger time-frequency region than the one shown here.}
\end{figure}

For this experiment, we set $q = 1$ and calculate the total energy on $R^q_m(C, f)$ for all $m = 1,2,3,4$ and all representations $C$ given by SST-STFT (which we also call \emph{first-order SST-STFT}), second-order SST-STFT, RM, and SST-QSTFT.  For SST-QSTFT we use Algorithm \ref{algorithm-2} with the parameters used for Figure 3. For all the transforms, we reassign the magnitude-squared of the STFT or QSTFT coefficients, in order to preserve the signal energy (see Remark \ref{remark-sstqstftdmag}). For first-order SST-STFT, second-order SST-STFT, and RM, we set the window $g = h^0$. For all transforms, $H=250$, $K=16384$, and $\gamma =0$ as before. Table \ref{table-second-order-sst-vs-sst-qstft} demonstrates that SST-QSTFT achieves superior ridge concentration as the SNR decreases, especially for the IFs $\phi'_1$ and $\phi'_4$ with higher absolute slope.  This underlies the point that the usage of chirped windows effectively exposes the presence of chirped signal content in heavy noise (as shown in Figure \ref{figure-four-transforms}).  By contrast, the second-order accurate reassignment frequency formula used for second-order SST-STFT is not as effective in detecting the local chirp rates $\{\sigma_m\}_{m=1}^4$ when noise is present.  Similarly, the reassignment time used in RM is not enough to compensate for the noise.  However, both second-order SST-STFT and RM improve over first-order SST-STFT, which is known to be ineffective when $|\sigma_m| \gg 0$\cite{oberlin2015second}.

\subsubsection{Application to analysis of gibbon calls}
In this section, we describe the potential of using SST-QSTFT together with a chirped window family to analyze audio recordings of animal calls.  In particular, we analyze a dataset of recordings of female Bornean gibbon great calls, recorded at the Stability of Altered Forest Ecosystems site in Sabah, Malaysia\cite{clink2017investigating}.  One application of interest is distinguishing the individual gibbons from each other based on features extracted from their calls.  

In what follows, we compute the SST-QSTFT of two calls from different gibbons, both with sampling rate $f_s = 44100$ Hz. To analyze these signals, we use the chirped window collection $\{h^w\}$ where $h^w(t) := h^0(t)\e^{\twopi \im \sigma^w t^2/2}$, $h^0$ and $H$ are as before, $K = 3L = 12288$, and where we select $\{\sigma^w\}_{w=1}^5 = \{ 50, 100, 200, 400, 600 \} \text{Hz}^2$ by visual inspection of a single STFT of the signal.  We use Algorithm \ref{algorithm-2} for window selection, with $\alpha = 0.5$, $A = B = 24$, $\Tstep = \Ystep = 4$, $\Tshift = \Yshift = 2$, and $P_T = P_Y = \{-1, 1\}$.

\begin{figure}
\includegraphics[width=18cm, height=6cm]{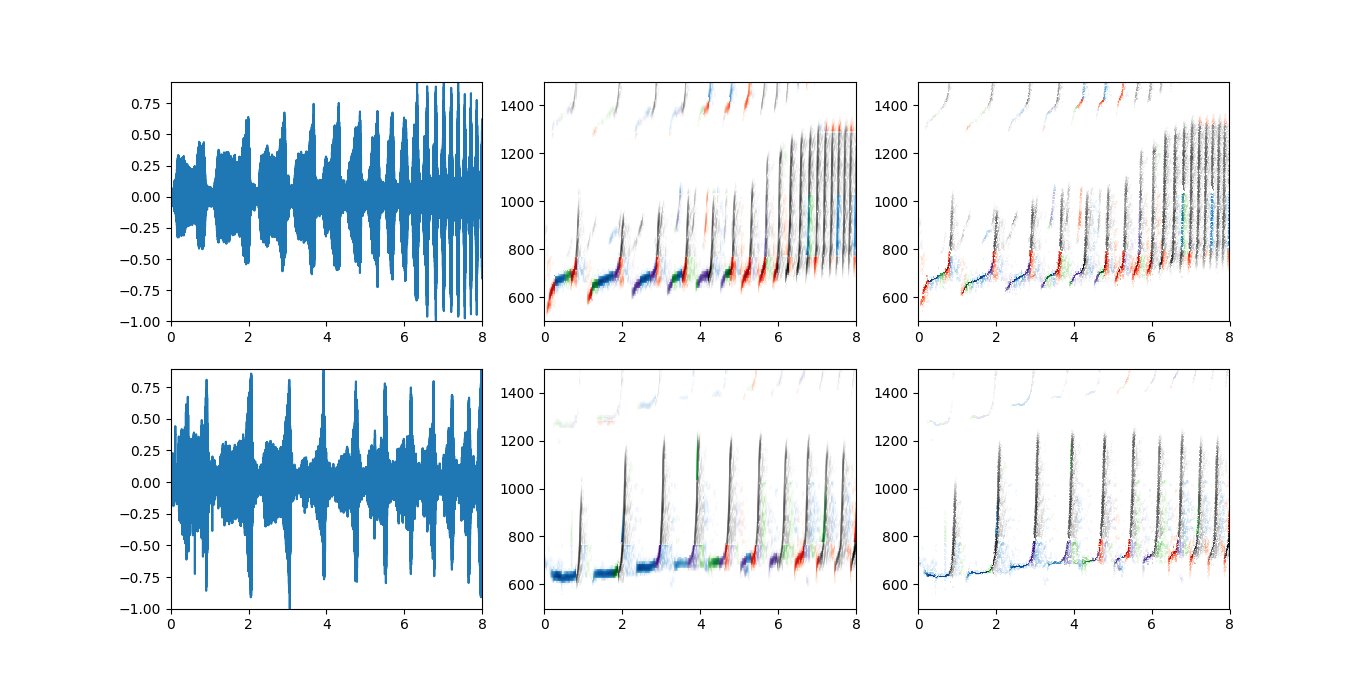}
\caption{\label{figure-gibbons} \tbf{Top row:} First gibbon. \tbf{Bottom row:} Second gibbon. \tbf{From left to right:} Original signal, QSTFT, and SST-QSTFT of the respective gibbon call.  Here, we restrict the frequency range of the plots between 500 and 1500 Hz, where the fundamental frequency information is captured.  Different colors correspond to the window chosen in each supertile for the calculation of the QSTFT.  Blue: $h^1$. Green: $h^2$. Purple: $h^3$. Red: $h^4$. Black: $h^5$.}
\end{figure}

We plot the QSTFT and SST-QSTFT of both signals in Figure \ref{figure-gibbons}.  Here, we see that the SST-QSTFT is a more concentrated time-frequency representation of the signal than the QSTFT.  Moreover, we may see also that the different curvature of the IF curves for the different gibbons is well-represented by the difference in optimal chirp rates.  In future work, we seek to use the QSTFT and SST-QSTFT to extract features from the gibbon calls that may be combined with a machine learning algorithm to distinguish the gibbons from one another.

\section{Conclusion}
We have developed a synchrosqueezing transform in the context of a quilted Gabor framework, enabling for improved adaptivity to the signal, and coming closer to the goal of an ideal time-frequency representation.  The SST-QSTFT yields an improved visualization of the IF information of the signal with adaptation to the time-frequency content, permitting for the precise isolation of diverse time-frequency events.  We have furthermore implemented a new algorithm for automatically and adaptively selecting optimal windows depending on the time-frequency content.  Reconstruction of each component is possible, even in the considered case of time-limited windows, by passing back to the original QSTFT coefficients.  Theoretical results demonstrate the concentration of the QSTFT around each IF curve and the closeness of the reassignment frequency to each true IF in both the continuous and discrete frameworks, as well as the accuracy of modes reconstruction in the continuous case.  Our numerical results show the effectiveness of the SST-QSTFT in adapting to the signal to achieve improved time-frequency concentration.

\acknowledgments

The authors were supported in part by ONR grant {N00014-16-1-2255}, as well as NSF grants {DMS-1418779} and {IIS-1631329}, and the GAANN fellowship grant {P200A120162}.  We thank Dena Jane Clink and Mark Grote of the UC Davis Department of Anthropology for providing us with the gibbon call dataset.  We also thank Dmytro Iatsenko for granting us permission to freely distribute the part of our Python code suite that we based on his MATLAB codes\footnote{Dmytro Iatsenko's codes can be found at \url{http://www.physics.lancs.ac.uk/research/nbmphysics/diats/tfr/}.} for SST\cite{iatsenko2015linear, iatsenko2013extraction}.  We used this Python code suite, including the part based on Dmytro Iatsenko's codes, to generate the results in this paper.

\bibliography{bibliography}
\bibliographystyle{spiebib}

\end{document}